\documentclass[a4paper, 11pt]{article}
\usepackage{stmaryrd}

\usepackage{amscd}
\usepackage{amssymb,amsmath,amsthm}
\usepackage{mathptmx}
\usepackage{mathrsfs}
\usepackage{color}
\usepackage{xspace}
\usepackage{bussproofs}

\usepackage{authblk}

\usepackage{cleveref}

\usepackage{enumerate}

\usepackage{relsize}

\DeclareSymbolFont{letters}{OML}{cmm}{m}{it}

\DeclareMathAlphabet{\mathcal}{OMS}{cmsy}{m}{n}

\usepackage{float}

\usepackage{tikz}
\usetikzlibrary{arrows}

\usetikzlibrary{decorations.markings}

\tikzstyle{index on}=[inner sep=2pt, white, circle, fill=black]
\tikzstyle{index off}=[inner sep=2pt, black, circle, draw]
\tikzstyle{index gray}=[inner sep=2pt, black, circle, fill=red!20]
\tikzstyle{opaque}=[fill=red,fill opacity=.1]

\newtheorem{theorem}{Theorem}[section]
\newtheorem{corollary}[theorem]{Corollary}
\newtheorem{lemma}[theorem]{Lemma}
\newtheorem{proposition}[theorem]{Proposition}

\newtheorem{example}[theorem]{Example}
\newtheorem{ctexample}[theorem]{Counter-example}

\newtheorem{fact}[theorem]{Fact}

\theoremstyle{definition}
\newtheorem{definition}[theorem]{Definition}

\usepackage{caption}

\newcommand{\Prop}{\ensuremath{\mathsf{Prop}}\xspace}
\newcommand{\MM}{\ensuremath{\mathfrak{M}}\xspace}

\newcommand{\IPC}{\ensuremath{\mathsf{IPC}}\xspace}
\newcommand{\CPC}{\ensuremath{\mathsf{CPC}}\xspace}

\newcommand{\IK}{\ensuremath{\mathsf{IK}}\xspace}
\newcommand{\IKt}{\ensuremath{\mathsf{IK}_\otimes}\xspace}

\newcommand{\IPL}{\ensuremath{\mathsf{IPL}}\xspace}
\newcommand{\TIPL}{\ensuremath{\mathsf{t}}\xspace}

\newcommand{\TGIPC}{\ensuremath{\mathsf{gtIPC}}\xspace}

\newcommand{\LL}{\ensuremath{\mathsf{L}}\xspace}
\newcommand{\tl}{\ensuremath{\mathsf{t}}\xspace}

\newcommand{\TIPC}{\ensuremath{\mathsf{tIPC}}\xspace}

\newcommand{\dblsetminus}{\mathbin{{\setminus}\mspace{-5mu}{\setminus}}}
\newcommand{\vvee}{\raisebox{1pt}{\ensuremath{\,\mathop{\mathsmaller{\mathsmaller{\dblsetminus\hspace{-0.23ex}/}}}\,}}}
\newcommand{\bigvvee}{\ensuremath{\mathop{\mathlarger{\dblsetminus}\hspace{-0.3ex}\raisebox{0.4pt}{$\mathlarger{\mathlarger{/}}$}}}\xspace}
\newcommand{\bigCup}{{\ensuremath{\mathlarger{\mathlarger{\mathlarger{\raisebox{-1.8pt}{$\Cup$}}}}}}\xspace}
\newcommand{\bigsCup}{{\ensuremath{\mathlarger{\mathlarger{\raisebox{-1pt}{$\Cup$}}}}}\xspace}

\title{Intermediate logics in the setting of team semantics}

\author[1]{Nick Bezhanishvili}
\author[2]{Fan Yang
}
\affil[1]{Institute for Logic, Language and Computation, University of Amsterdam, The Netherlands}
\affil[2]{Department of Mathematics and Statistics, University of Helsinki, Finland}

\date{}

\begin{document}

\maketitle


\begin{center}
{\small \em Dedicated to Dick de Jongh on the occasion of his 80th birthday.}
\end{center}

\vspace{\baselineskip}

\begin{abstract}
Several authors have recently defined intuitionistic logic based on team semantics (\TIPC).  In this paper we provide two alternative approaches to  intermediate logics in the team semantics setting. We do this by modifying \TIPC with  axioms written with 
two different versions of disjunction in the logic, a local one and global one. We prove a characterization theorem in the first approach
and we introduce a generalized team semantics in the second one.
\end{abstract}

\section{Introduction}



One of the important themes in the research of Dick de Jongh has been the study of intuitionistic and other  intermediate logics. De Jongh  made a number of important contributions in this field  and is, in fact, one of the founders of this area
(together with 
Gabbay, Troelstra, Hosoi, Jankov, Kuznetsov and others). Together with Gabbay, de Jongh defined an important class of intermediate logics 
nowadays called the \emph{Gabbay-de Jongh logics} \cite{G-deJ74}. With Troelstra he developed what is currently known as the \emph{discrete duality for Heyting 
algebras} \cite{T-deJ66}.  
In  1960's, de Jongh proved in his PhD thesis  
 possibly the first De Jongh Theorem, 
the characterization of intuitionistic propositional calculus $\mathsf{IPC}$ among other intermediate logics 
in terms of the 
\emph{Kleene Slash} \cite{deJonghPhD} (see also \cite{Bezh04}).  
In order to obtain this characterization de Jongh introduced  a type of formulas that we nowadays  call \emph{de Jongh formulas}. This is based on special ``coloring'' of Kripke structures, which was later developed by several authors to a construction known as the $n$-universal models for intuitionistic logic (see, e.g.,  \cite{ZhaCha_ml} and \cite{BezhanishviliPhD}).
Independently, using algebraic techniques, Jankov \cite{Jan63} introduced a method of, what we now call, \emph{Jankov} or \emph{splitting formulas} for Heyting algebras. This allowed Jankov to construct intermediate logics without the finite model property  and to prove that the lattice of intermediate logics has the cardinality that of the continuum \cite{Jan63,Jan68}.  It turnes out that Jankov and de Jongh formulas are two sides of the same coin. In particular, the Jankov formula of a finite subdirectly irreducible 
Heyting algebra $A$ is semantically equivalent to the de Jongh formula of the finite rooted Kripke frame 
dual to $A$ (for details of this connection we refer to \cite[Section 3.2.2]{BezhanishviliPhD}). Because of this, these formulas
are also known as \emph{Jankov-de Jongh formulas} \cite[Remark 3.3.5]{BezhanishviliPhD}. An alternative frame-theoretic description of such formulas can be found in \cite[Theorem 9.39]{ZhaCha_ml}. 
From the modern perspective, de Jongh and Jankov formulas axiomatize splittings and their joins in the lattice
of intermediate logics. 
The analogues of de Jongh formulas for transitive modal logics were developed by Fine \cite{Fin74} and the analogues of Jankov formulas for $\mathsf{K4}$-modal algebras by Rautenberg \cite{Rau80}. For further generalizations of such kind of formulas and their impact on the theory of modal and intermediate  logics we refer to \cite[Chapter 9]{ZhaCha_ml} and \cite{BezhJank}. 
De Jongh formulas and the method of universal models has become a standard tool in the study of  various fragments of intermediate  and other non-classical logics \cite{BJTZ15, BCJ17}. 
As a tribute to Dick de Jongh on the occasion of his 85th birthday, in the first part of this paper we  give an overview of the method of de Jongh formulas and sketch some of their new applications.  

In the second, main part of the paper, we study intermediate logics in a new domain of team semantics. We present two different approaches. In the first and simpler approach, we provide yet another new application for de Jongh formulas. The second approach is our first attempt to introduce a generalized team semantics in which intermediate logics can be naturally defined. Thus we lift intermediate logics (de Jongh's paradise) to a new setting. 
We hope that this paper demonstrates again how the ideas developed by Dick de Jongh decades ago could influence new research lines to this day. We now discuss our motivation and contributions in detail.

Team semantics was originally introduced by
Hodges \cite{Hodges1997a,Hodges1997b} to characterize notions of dependence  in Hintikka and Sandu's {\em independence friendly-logic} \cite{Hintikka98book,HintikkaSandu1989}. This framework was later developed by V\"{a}\"{a}n\"{a}nen in his {\em dependence logic} \cite{Van07dl}. Independently, Ciardelli and Roelofsen \cite{InquiLog} developed {\em inquisitive logic}, which turns out to be a logic also based on team semantics\footnote{The fact that inquisitive logic essentially also adopts team semantics was first observed by Dick de Jongh and Tadeusz Litak in a private conversation with the second author in 2011.} 
 and was shown \cite{Ciardelli2015,Yang_dissertation} to be actually a variant of propositional dependence logic \cite{VY_PD}.  Both dependence and inquisitive logic can be understood as conservative extensions of classical logic.  Recently several authors have defined different intuitionistic logic-based dependence and inquisitive logic \cite{CiardelliIemhoffYang20,Holliday2020,Puncochar16,Puncochar17}. 
The starting point of this paper is the team-based intuitionistic logic of \cite{CiardelliIemhoffYang20}. 

The key idea of team semantics is that formulas are evaluated in a model with respect to {\em teams}, which, in the intuitionistic Kripke semantics setting, are defined as {\em sets} $t\subseteq W$ of possible worlds of an intuitionistic Kripke model $\mathfrak{M}=(W,R,V)$. Such sets $t$ correspond to multitudes of state of affairs, or they express certain  uncertainly of the current state of affairs.
We also extend the standard language 
of intuitionistic propositional logic (\textsf{IPC}) 
with  a natural disjunction  $\vvee$ (called {\em global disjunction}) 
on the team level (i.e., the level of $\wp(W)$). The disjunction $\vee$ of \textsf{IPC} is thus referred to as {\em local disjunction}. {We denote our full language as $\mathop{[\bot,\wedge,\vee,\vvee,\to]}$, and the standard language of \IPC as $[\bot,\wedge,\vee,\to]$. Formulas in the standard language $[\bot,\wedge,\vee,\to]$ (i.e., $\vvee$-free formulas) are referred to as {\em standard formulas}.}
The logical system  (denoted by \textsf{tIPC}) as defined in \cite{CiardelliIemhoffYang20} is a conservative extension of \IPC for standard formulas with additional axioms, including \IPC axioms  for the language {with connectives $\bot,\wedge,\vvee,\to$}, and the \textsf{Split} axiom
\[\alpha\to(\phi\vvee\psi)\to(\alpha\to \phi)\vvee(\alpha\to\psi),\tag{\textsf{Split}}\]
where $\alpha$ is an arbitrary standard formula  (thus \textsf{tIPC} is not closed under uniform substitution).



We provide two alternative approaches to  intermediate logics in this setting, by modifying \TIPC with  axioms written with  two different disjunctions. 
Given an intermediate logic $ \textsf{L}=\textsf{IPC}\oplus\Delta$ with $\Delta$ a set of axioms in the standard language, the first approach defines an intermediate team logic $\textsf{tL}$ by adding to \TIPC all instances of $\Delta$ with propositional variables substituted with standard formulas only, and closing the resulting set under Modus Ponens.
This definition has appeared also in a recent work \cite{Quadrellaro21}. A similar but subtly different definition was given in \cite{Puncochar21}. We show, by using a disjunctive normal form $\bigvvee_{i\in I}\alpha_i$ (with each $\alpha_i$ a standard formula) of the logic, that $n$-universal models of \IPC still behave as universal models for \TIPC, and thus de Jongh formulas for team-based intuitionistic logic can be defined in the usual manner. Moreover, a large class of intermediate axioms (including de Jongh formulas) for \IPC still characterize the same class of frames in the team semantics setting as in the standard (single-world semantics) setting, in the sense that  if an intermediate logic \textsf{L} (over the standard language) is complete with respect to a class $\textsf{F}$ of frames, then \textsf{tL} is also complete with respect to $\textsf{F}$, provided that $\textsf{L}$ has disjunction property or $\textsf{L}$  is canonical. 

In the second approach, we seek appropriate ways to obtain logics in which  
the \textsf{Split} axiom of \textsf{tIPC} can be replaced with other (weaker) axioms (written in the full language 
$[\bot,\wedge,\vee,\vvee,\to]$). 
It is fair to say that the \textsf{Split} axiom is the hallmark of team-based propositional logics (as, e.g., it is a crucial axiom for establishing the important disjunctive normal form $\bigvvee_{i\in I}\alpha_i$ for these logics). Weakening the \textsf{Split} axiom and at the same time keeping certain unique features of team semantics in some meaningful manner is thus not a trivial task. In this paper, we make a first attempt in this direction.
 We start by exploring a connection between \TIPC and intuionistic modal logic \IK \cite{WolterZakhIML99}, given by the observation that every (intuitionistic) Kripke model $\mathfrak{M}=(W,R,V)$ of \TIPC gives rise to a (full) powerset model $\mathfrak{M}^\bullet=(\wp(W),\supseteq, R^\circ,\cup,\varnothing,V^\circ)$, where $R^\circ$ and $V^\circ$ are liftings of $R$ and $V$, and \TIPC can thus be given a single-world semantics interpretation over $\mathfrak{M}^\bullet$. In this setting the local disjunction $\vee$ is naturally understood as a binary diamond modality. Such a connection was studied also in the modal dependence logic setting in \cite{Yang17MD}. Our powerset model construction also corresponds to the set-lifting as discussed in \cite{vanBenthem21}, and our local disjunction $\vee$ is essentially  the binary diamond $\langle \sup\rangle$ defined in \cite{vanBenthem21}.
Inspired by these connections, we then introduce the so-called {\em general intuitionistic team Kripke  models} (which generalizes the full powerset models). Over these models we define a generalized team intuitionsitic Kripke semantics in which the \textsf{Split} axiom is no longer sound.
The \textsf{Split} axiom essentially enforces the underlying structure of teams (i.e., the part of the full powerset model $\mathfrak{M}^\bullet$ characterizing teams) to be $(\wp(W),\supseteq,\cup,\varnothing)$. In our approach, weakening the  \textsf{Split}  axiom  amounts to changing this structure to $(\wp(W),\succcurlyeq,\Cup,\varnothing)$ with a bounded semilattice $(\wp(W),\Cup,\varnothing)$ and the associated partial order $\preccurlyeq$ satisfying certain constraints,
and thus generalizing the standard team semantics. 
In  related work Pun\u{c}och\'{a}\u{r}  \cite{Puncochar16} generalized  the semantics by essentially considering an underlying structure $(I,\supseteq,\cup,\varnothing)$ with a set $I\subseteq \wp(W)$ of teams equipped with certain topology as the underlying structure for teams, resulting in different intermediate team logics from ours (in particular, these logics do not have the local  and global disjunction at the same time). Further Pun\u{c}och\'{a}\u{r} \cite{Puncochar17} considered a more abstract generalization of team semantics, by treating teams as primitive entities in arbitrary bounded semilattices $(A,\Cup,0)$. The resulting logic does satisfy the \textsf{Split} axiom in case the semantics is persistent.
We will provide more comparisons of our work with  \cite{Puncochar16,Puncochar17} in the concuding section, Section 5. 
There are also similarities between our framework and the framework of possibility semantics \cite{Holl-Pos}. In particular, the disjunction of possibility semantics over canonical possibility models coincides with the local disjunction in our setting, as shown in \cite[Corollary~5.5]{BH20}. 

Under the generalized team semantics,  two \TIPC axioms of the local disjunction $\vee$, the monotonicity and (weak) elimination axiom, turn out  to be not sound either. 
We also study frame conditions under which the soundness of these two axioms as well as the \textsf{Split} axiom can be recovered. 
We show that distributivity  is a sufficient condition for validating all these axioms, and over finite classical frames, distributivity is also a necessary condition for the two axioms of the local disjunction. We leave it  for future work to identify in the general case the necessary and sufficient conditions for validating these axioms. We hope that the  conditions provided in this paper already give  some new insights on the peculiar 
properties of  team semantics and of the logic \TIPC.




This paper is organized as follows. In Section 2, we provide an account of intuitionistic propositional logic \IPC, universal models 
and de Jongh formulas in the standard (single-world Kripke semantics) setting. In Section 3, we recall from \cite{CiardelliIemhoffYang20} intuitionistic propositional logic \TIPC based on team semantics, and investigate intermediate team logics through the first approach discussed above. Section 4 explores the second approach.  We end in Section 5 by some concluding remarks and open problems.

\section{Intuitionistic logic and de Jongh formulas}

In this section, we recall briefly intuitionistic logic and de Jongh formulas in the standard setting. For more detailed discussions, the reader is referred to, e.g., \cite{BezhanishviliPhD,ZhaCha_ml}.

Fix a set $\Prop$ of propositional variables. The language of {\em intuitionistic propositional logic} 
is defined as
\[\phi::=p\mid \bot\mid\phi\wedge\phi\mid \phi\vee\phi\mid\phi\to\phi.\]
We write $\neg\phi:=\phi\to\bot$. 
We consider the usual {\em intuitionistic propositional calculus} (\IPC), which is classical propositional calculus (\CPC) without the law of excluded middle ($\phi\vee\neg\phi$) or the double negation elimination axiom ($\neg\neg\phi\to\phi$).  We write $\vdash_{\IPC}\phi$ or simply $\vdash\phi$ if $\phi$ is a theorem.

An {\em (intuitionistic) Kripke model} is a tuple $\mathfrak{M}=(W,R,V)$ such that $W$ is a non-empty set of {\em worlds} (or {\em nodes} or {\em points}), $R\subseteq W\times W$ is a partial order, and $V:\Prop\to \wp(W)$ is a {\em persistent} valuation, i.e., $w\in V(p)$ and $wRv$ imply $v\in V(p)$. The underlying pair $\mathfrak{F}=(W,R)$ is called an {\em (intuitionistic) Kripke frame}.  If $wRv$, then $v$ is called a {\em successor} of $w$.  An {\em immediate successor} $v$ of $w$ is a proper successor of $w$ such that there is no other point $u$ with $wRu$ and $uRv$. Denote by $R\lceil w\rceil$ (or simply $\lceil w\rceil$) the set of all immediate successors of $w$. 
A point $w$ is called an {\em endpoint} if $\lceil w\rceil=\varnothing$. Generated subframes or submodels, and p-morphisms between models or frames are defined as usual. We write $\mathfrak{M}_w$ for the submodel of $\mathfrak{M}$ generated by a point $w$ in $\mathfrak{M}$. The satisfaction relation $\mathfrak{M},w\models_{\IPL}\phi$ (or simply $\mathfrak{M},w\models\phi$) is defined inductively in a standard way; in particular, 
\begin{itemize}
\item $\mathfrak{M},w\not\models \bot$ ~~always;
\item $\mathfrak{M},w\models \phi\to\psi$ ~~iff~~ for any $v\in W$ with $wRv$, $\mathfrak{M},v\models \phi$ implies $\mathfrak{M},v\models \psi$.
\end{itemize}
The {\em persistence} property generalizes easily to arbitrary formulas; that is, $\mathfrak{M},w\models\phi$ and $wRv$ imply $\mathfrak{M},v\models\phi$ for  arbitrary formulas $\phi$.

Given $n$ propositional variables $p_1,\dots,p_n\in\Prop$. A model $\mathfrak{M}=(W,R,V)$ with the valuation $V$ restricted to $p_1,\dots,p_n$ is referred to as an $n$-model. It is often convenient to describe the valuation  in an $n$-model by using {\em colors}. The {\em ($n$-)color} of a point $w\in W$ is defined as a $0$-$1$-sequence $col(w)=c_1\dots c_n$ such that $c_i=1$ iff $w\models p_i$, and $c_i=0$ iff $w\not\models p_i$. We write $c_1\dots c_n\leq c_1'\dots c_n'$ if $c_i\leq c_i'$ for each $i$; and write $c_1\dots c_n<c_1'\dots c_n'$ if $c_1\dots c_n\leq c_1'\dots c_n'$ and $c_1\dots c_n\neq c_1'\dots c_n'$. 

De Jongh formulas can be defined over universal models of \IPC.
For any natural number $n\in \mathbb{N}$,  the $n$-universal model $\mathcal{U}(n)=(W,R,V)$ of \IPC is defined inductively in layers as follows:
\begin{itemize}
\item The first layer consists of $2^n$ nodes of distinct $n$-color.
\item Suppose the $\leq k$ layers have been already defined. We define the ($k+1$)th layer as follows:
\begin{itemize}
\item For every node $w$ in the $k$th layer, and every color $c<col(w)$, add a new node $v$ such that $col(v)=c$ and $\lceil v\rceil =\{w\}$.
\item For every finite set $A$ of pairwise $R$-incomparable nodes in layers $\leq k$ with at least one element from the $k$th layer, and for every $n$-color $c\leq col(w)$ for all $w\in A$, add a new node $v$ such that $col(v)=c$ and $\lceil v \rceil=A$.
\end{itemize}
\end{itemize}

Every finite rooted $n$-model can be found as a p-morphic image of a generated submodel of the $n$-universal model, as the next lemma shows. 
We include here also a detailed proof of this well-known fact, proof  of which however does not seem to have been well documented in the literature. 

\begin{lemma}\label{Un_universality_lm}\
\begin{enumerate}[(i)]
\item Every finite rooted Kripke $n$-model $\mathfrak{M}$ can be mapped p-morphically onto a unique generated submodel of $\mathcal{U}(n)$.
\item Every finite rooted Kripke frame $\mathfrak{F}$ is isomorphic to a generated subframe of $\mathcal{U}(n)$ for some $n\leq |\mathfrak{F}|$.
\end{enumerate}
\end{lemma}
\begin{proof}
(i) We prove the claim by induction on the depth of $\mathfrak{M}$ with root $r$. If $\mathfrak{M}$ is a singleton $\{r\}$, then obviously $r$ can be mapped p-morphically to a unique endpoint $w$ in $\mathcal{U}(n)$ of the same color. Now, suppose $\lceil r\rceil\neq \varnothing$. By induction hypothesis, each $v\in \lceil r\rceil$ can be mapped, via some p-morphism $f_v$ onto a unique submodel of $\mathcal{U}(n)$ generated by some point $w_v$. Let $A$ be the set of minimal points in $\{w_v\mid v\in\lceil r\rceil\}$.

If $A$ is a singleton $\{w_{v_0}\}$, then we distinguish two cases. If $col(r)=col(w_{v_0})$, then it is easy to see that the function $f=f_{v_0}\cup\{(r,w_{v_0})\}$ is a (unique) onto p-morphism. If $col(r)<col(w_{v_0})$, then by the construction of $\mathcal{U}(n)$, there is a (unique) point $w_r$ in $\mathcal{U}(n)$ with $col(w_r)=col(r)$ and $\lceil w_r\rceil=\{w_{v_0}\}$. It is easy to see that the function $f=f_{v_0}\cup\{(r,w_r)\}$ is a (unique) onto p-morphism.

If $A$ is a set of pairwise incomparable nodes, then by the construction of $\mathcal{U}(n)$, there is a (unique) point $w_r$ in $\mathcal{U}(n)$ with $col(w_r)=col(r)$ and $\lceil w_r\rceil=A$. Thus, the function $f=\{(r,w_r)\}\cup\bigcup_{w_v\in A}f_v$ is a p-morphism. We also have $f_v\upharpoonright \mathcal{U}(n)_u=f_u$ for every $w_v\in A$ and $u\in \lceil r\rceil$ with $w_vRw_u$. Thus, it is not hard to see that $f$ is the unique onto p-morphism.



(ii) Let  $\mathfrak{F}=(W,R)$ be a finite rooted frame. Introduce a new propositional variable $p_w$ for each point $w$ in $W$, and define a valuation $V$ on $\mathfrak{F}$ by setting $V(p_w)=R(w)$. Putting $n=|W|$, by item (i), $(\mathfrak{F},V_{\mathfrak{F}})$ can be mapped p-morphically via $f$ onto a  (unique) generated submodel of $\mathcal{U}(n)$ with root $w$. Since every point in $(\mathfrak{F},V_{\mathfrak{F}})$ has a distinct color, the p-morphism $f$ is an isomorphism, namely that $(\mathfrak{F},V_{\mathfrak{F}})$ is isomorphic to the generated submodel $\mathcal{U}(n)_w$. 
Note that the underlying frame of $\mathcal{U}(n)_w$ may also be a generated subframe of some $\mathcal{U}(m)$ with $m<n$.
\end{proof}

It follows immediately from the above lemma that universal models are  the counter-models for all non-theorems  of \IPC, as (finite  rooted) counter-models of every non-theorem of \IPC can be mapped p-morphically onto the universal model. 

\begin{theorem}\label{Un_prop}
For every $n$-formula $\phi$, we have that
\(\mathcal{U}(n)\models\phi\iff \vdash_{\IPC}\phi.\)
\end{theorem}

We are now ready to recall the definition of de Jongh formulas.

\begin{definition}[De Jongh formulas] 
Let $w$ be a point in $\mathcal{U}(n)$. Define $n$-formulas $\phi_w$ and $\psi_w$ by induction on the depth of $w$ as follows:
\begin{itemize}
\item If $w$ is an endpoint, then define
\[\phi_w:=\bigwedge_{w\models p_i}p_i\wedge\bigwedge_{w\not\models p_i}\neg p_i\quad\text{ and }\quad\psi_w:=\neg \phi_w.\]
\item If $w$ is not an endpoint (i.e., $\lceil w\rceil\neq\varnothing$), then define
\begin{align*}
\phi_w&:=\bigwedge_{w\models p_i} p_i\wedge \big(\bigvee \textsf{np}(w)\vee\bigvee_{v\in \lceil w\rceil}\psi_v\to \bigvee_{v\in \lceil w\rceil}\phi_v\big)\\
\text{ and }\quad\psi_w&:=\phi_w\to\bigvee_{v\in \lceil w\rceil}\phi_v,
\end{align*}
where
\(\textsf{np}(w):=\{p_i\mid w\not\models p_i\text{ and }v\models p_i\text{ for all }v\in \lceil w\rceil,~1\leq i\leq n\}.\)
\end{itemize}

\end{definition}

The following theorem describes the important property of de Jongh formulas.

\begin{theorem}\label{de_jong_form_prop}
For every points $w,u$ in $\mathcal{U}(n)$,
\begin{enumerate}[(i)]
\item $\mathcal{U}(n),u\models\phi_w\iff wRu$;
\item $\mathcal{U}(n),u\not\models\psi_w\iff uRw$.
\end{enumerate}
\end{theorem}
\begin{proof}
Item (ii) is a consequence of item (i), since
\begin{align*}
u\not\models\psi_w&\iff \exists x: uRx,~x\models\phi_w\text{ and }\forall v\in\lceil w\rceil: x\not\models\phi_v\tag{in case $w$ is an endpoint, the second part holds trivially}\\
&\iff  \exists x: uRx,~wRx\text{ and }\forall v\in \lceil w\rceil: \neg vRx\tag{by item (i)}\\
&\iff \exists x: uRx\text{ and }w=x,\text{ i.e., }uRw\tag{since $R$ is a partial order}
\end{align*}

We now prove item (i) by induction on the depth of $w$. If $w$ is an endpoint (i.e., $w$ is  in the first layer of $\mathcal{U}(n)$), the direction ``$\Longleftarrow$" is trivial; for the direction ``$\Longrightarrow$", it is easy to see that $\mathcal{U}(n),u\models_{\IPL}\phi_w$ implies that $u$ and all its successors have the same color as $w$, which can only happen when $wRu$ by the construction of $\mathcal{U}(n)$. 

Now, suppose $w$ has proper successors. We first show the direction ``$\Longleftarrow$". Suppose $wRu$. Clearly, $u\models\bigwedge_{w\models p_i}p_i$. If $u=w$, then by induction hypothesis, we have that $w\not\models\bigvee_{v\in \lceil w\rceil}\psi_v$. If $u$ is a proper successor  of $w$, then $u$ is a successor of some point in $\lceil w\rceil$, which by induction hypothesis implies that $u\models \bigvee_{v\in \lceil w\rceil}\phi_v$. Putting everything together entail  that $u\models\phi_w$.

For the converse direction ``$\Longrightarrow$", assume  that $\neg (wRu)$. We will show that $u\not\models\phi_w$. Suppose also that $u\models p_i$ for all $p_i$ with $w\models p_i$, meaning that $col(u)\geq col(w)$. 
If there exists a successor $x$ of $u$  such that $\neg (wRx)$ and  $\neg (xRv_0)$ for some $v_0\in \lceil w\rceil$, then by induction hypothesis, we have that  $x\not\models \bigvee_{v\in \lceil w\rceil}\phi_v$ (for otherwise $wRvRx$ for some $v\in \lceil w\rceil$) and $x\models \bigvee_{v\in \lceil w\rceil}\psi_v$. Thus, $u\not\models\phi_w$. 

Now, suppose that no such point exists, namely that 
\[\forall x\in R[u]:\text{ either }wRx\text{ (i.e., $w\in R[\lceil w\rceil]$)}\text{ or } \lceil w\rceil\subseteq R[x].\]
That is, every successor $x$ of $u$ can see all points in $\lceil w\rceil$ in the future, or have seen a point in $\lceil w\rceil$ in the past. By the construction of $\mathcal{U}(n)$, this implies that there exist successors $x$ of $u$ such that  $\lceil x\rceil=\lceil w\rceil$. Take such an $x$ with $x\neq w$. Such an $x\neq w$ must exist, for otherwise, by the construction of $\mathcal{U}(n)$, we must have that $\lceil u\rceil$ is a singleton of a predecessor of $x$ and $col(u)<col(w)$, which is a contradiction.

Now, again by the construction of $\mathcal{U}(n)$, we must have that $col(x)\neq col(w)$. Since $uRx$, we also have $col(x)\geq col(u)\geq col(w)$. Thus $col(x)>col(w)$, meaning that there is some $p_j$ such that $w\not\models p_j$ but $x\models p_j$. For all $v\in \lceil w\rceil$, since $xRv$, we have that $v\models p_j$. It then follows that $p_j\in \textsf{np}(w)$ and $x\models \bigvee\textsf{np}(w)$. On the other hand,  by induction hypothesis, we have that $x\not\models\bigvee_{v\in \lceil w\rceil}\phi_v$. Hence, $u\not\models\phi_w$.
\end{proof}

We can now prove the Jankov-de Jongh theorem using de Jongh formulas and universal models. Note that Jankov's proof \cite{Jan63} was purely algebraic 
(we refer to  \cite{BezhJank} for more details on Jankov's approach). Yet another perspective to this formulas  when one introduces new variable for each element of the frame 
is by Fine \cite{Fin74} in the modal case and by Zakharyaschev in the intuitionistic one (see \cite[Chapter 9]{ZhaCha_ml}).


\begin{theorem}[Jankov-de Jongh]
For any finite rooted frame $\mathfrak{F}$, there exists a formula $\chi(\mathfrak{F})$ such that  for every frame $\mathfrak{G}$,
\[\mathfrak{G}\not\models\chi(\mathfrak{F}) \iff \mathfrak{F} \text{ is a generated subframe of a p-morphic image }\mathfrak{G}.\]
\end{theorem}
\begin{proof}
By Lemma \ref{Un_universality_lm}, $\mathfrak{F}$ is isomorphic to a generated subframe of $\mathcal{U}(n)$ with root $w$ and $n\leq \mathfrak{F}$. Let $\chi(\mathfrak{F})=\psi_w$.
``$\Longleftarrow$": Suppose $\mathfrak{F}$ is a generated subframe of a p-morphic image of $\mathfrak{G}$. By Theorem \ref{de_jong_form_prop}, $\mathcal{U}(n),w\not\models\psi_w$. Thus $\mathfrak{F}\not\models\psi_w$, which further implies that $\mathfrak{G}\not\models\psi_w$.

``$\Longrightarrow$": Suppose $(\mathfrak{G},V)\not\models\psi_w$, where $(\mathfrak{G},V)$ is some $n$-model on $\mathfrak{G}$. By Lemma \ref{Un_universality_lm}(i), $(\mathfrak{G},V)$ can be mapped p-morphically onto a generated submodel of $\mathcal{U}(n)$ with root $u$. Thus, $\mathcal{U}(n),u\not\models\psi_w$. By Theorem \ref{de_jong_form_prop}, this means that $uRw$. Since $\mathfrak{F}$ is isomorphic to the  frame of $\mathcal{U}(n)_w$, it follows that $\mathfrak{F}$ is a generated subframe of the underlying frame of $\mathcal{U}(n)_u$, which is a p-morphic image of $\mathfrak{G}$.
\end{proof}

The algebraic counter-part of the next lemma is the \emph{congruence extension property}, which in particular implies that for a Heyting algebra $A$ we have ${\bf HS}(A) = {\bf SH}(A)$, where ${\bf H}$ and ${\bf S}$ are the operations of taking subalgebras and homomorphic images, respectively (see, e.g., \cite{Esa19}).

\begin{lemma}
For any frames $\mathfrak{F}$ and $\mathfrak{G}$, $\mathfrak{F}$ is a generated subframe of a p-morphic image of $\mathfrak{G}$, iff $\mathfrak{F}$ is a p-morphic image of a generated subframe of $\mathfrak{G}$.
\end{lemma}
\begin{proof}
If $\mathfrak{F}$ is a generated subframe of $f[\mathfrak{G}]$ for some p-morphism $f$, then $f^{-1}[\mathfrak{F}]$ is a generated subframe of $\mathfrak{G}$ and $f[f^{-1}[\mathfrak{F}]]=\mathfrak{F}$. Conversely, if $\mathfrak{G}_0=(W_0,R)$ is a generated subframe of $\mathfrak{G}=(W,R)$ and $f[\mathfrak{G}_0]=\mathfrak{F}$ for some p-morphism $f$, then consider the extension $\mathfrak{H}=(W',R')$ of $\mathfrak{F}=(f[W_0],R')$,  defined as
\begin{itemize}
\item $W'=f[W_0]\cup (W\setminus W_0)$;
\item for any $w,v\in W\setminus W_0$, $wRv$ iff $wR'v$;
\item for any $w\in W\setminus W_0$ and $v'\in f[W_0]$, $wR'v')$ iff $wRv$ for some $v\in W_0$ such that $f(v)=v'$.
\end{itemize}
Let $g:\mathfrak{G}\to\mathfrak{H}$ be an extension of $f$ that is the identity function on $W\setminus W_0$. It is easy to show that $g$ is an onto p-morphism.
\end{proof}

\begin{corollary}
For any finite rooted frame $\mathfrak{F}$, there exists a formula $\chi(\mathfrak{F})$ such that  for every frame $\mathfrak{G}$,
\[\mathfrak{G}\not\models\chi(\mathfrak{F}) \iff \mathfrak{F} \text{ is a p-morphic image of a  generated subframe of }\mathfrak{G}.\]
\end{corollary}

\begin{figure}[t]
\begin{center}
\begin{tikzpicture}[scale=1, transform shape]

  \node[draw, circle, scale=0.5, fill=black, label={[left]\small $w_0$}] (0) at (0,0) {};
  
  \node at (0.2,0.4) {$p$};

  \node[draw, circle, scale=0.5, fill=black, label={[right]\small $w_1$}] (1) at (1.5,0) {};
 
  \node[draw, circle, scale=0.5, fill=black, label={[left]\small $w_2$}] (2) at (0,-1) {};
   \node[draw, circle, scale=0.5, fill=black, label={[right]\small $w_3$}] (3) at (1.5,-1) {};

   \node[draw, circle, scale=0.5, fill=black] (4) at (0,-2) {};    
   
      \node[draw, circle, scale=0.5, fill=black] (5) at (1.5,-2) {};
         \node[draw, circle, scale=0.5, fill=black] (6) at (0,-3) {};
  \node[draw, circle, scale=0.5, fill=black] (7)  at (1.5,-3) {};
         
         
\draw[line width=.75pt] (0) -- (3);
\draw[line width=.75pt] (1) -- (4);
\draw[line width=.75pt] (2) -- (5);
\draw[line width=.75pt] (3) -- (6);
\draw[line width=.75pt] (4) -- (7);
\draw[line width=.75pt] (5) -- (0,-4);
\draw[line width=.75pt] (6) -- (1.5,-4);
\draw[line width=.75pt] (7) -- (0.75,-4);
\draw[line width=.75pt] (0,0) -- (0,-3.3);
\draw[line width=.75pt] (1.5,0) -- (1.5,-3.3);
\draw[dashed, line width=.75pt] (0,-3.3) -- (0,-3.8);
\draw[dashed, line width=.75pt] (1.5,-3.3) -- (1.5,-3.8);

  \end{tikzpicture}
  \caption{Rieger-Nishimura ladder}
\label{Rieger-Nishimura}
\end{center}
\end{figure}
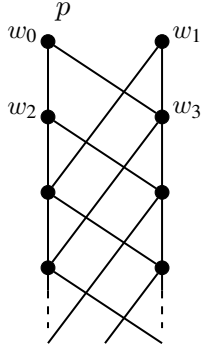

Given a set $\Delta$ of formulas, we write $\IPC\oplus\Delta$ for the set of formulas (or the logic) obtained by closing the set $\IPC\cup\Delta$ under uniform substitution and Modus Ponens. A logic $\LL$ is called an {\em intermediate logic} if $\IPC\subseteq \LL\subseteq \CPC$.  A large class of intermediate logics can be axiomatized by de Jongh formulas or sets of de Jongh formulas. These are exactly  the splititng and join splitting logics in the lattice of intermediate logics \cite[Section 10.5]{ZhaCha_ml}. 
We will now write explicitly some of the axiomatizations of the well-known intermediate logics in terms of de Jongh formulas. Recall that the 1-universal model of $\IPC$ is the Rieger-Nishimura ladder drawn in Figure \ref{Rieger-Nishimura} (see e.g., 
\cite[Section 8.7]{ZhaCha_ml} and \cite[Section 3.2]{BezhanishviliPhD}). We label it by the elements $w_i$ for $i\in \omega$. 
Recall also that the {\em Rieger-Nishimura polynomials} are given by the following
recursive definition: \index{Rieger-Nishimura! polynomials}
\begin{enumerate}
\item $g_0(p):= p$,
\item $g_1(p):= \neg\, p$,
\item $f_1(p):= p\vee \neg p$,
\item $g_2(p):= \neg\,\neg\, p$,
\item $g_3(p):= \neg\,\neg\, p\to p$,
\item $g_{n+4}(p):= g_{n+3}(p)\to (g_n(p)\vee g_{n+1}(p))$,
\item $f_{n+2}(p):= g_{n+2}(p)\vee g_{n+1}(p)$. \index{$g_k(p)$} 
\end{enumerate}
Then it is well known that the upsets of each $w_i$ are defined by the  Rieger-Nishumura poynomials. 
In particular, For every $k\in \omega$ we have:
\begin{itemize}
\item $R(w_k)=\{w\in \mathcal{U}(1):w\models g_k(p)\}$,

\item $R(w_k)\cup R(w_{k-1})=\{w\in \mathcal{U}(1):w\models f_k(p)\}$.
\end{itemize}

\begin{figure}[t]
\begin{center}
\begin{tikzpicture}[scale=1, transform shape]

  \node[draw, circle, scale=0.4, fill=black] (f0) at (0,0) {};

  \node[draw, circle, scale=0.4, fill=black] (f1) at (0,1) {};
  
 \node at (0,-0.5) {$\mathfrak{F}$}; 
         
\draw[line width=.75pt] (f0) -- (f1);

  \node[draw, circle, scale=0.4, fill=black] (g0) at (3,0) {};

  \node[draw, circle, scale=0.4, fill=black] (g1) at (2.3,1) {};
    \node[draw, circle, scale=0.4, fill=black] (g2) at (3.7,1) {};
  
 \node at (3,-0.5) {$\mathfrak{G}$}; 
         
\draw[line width=.75pt] (g0) -- (g1);
\draw[line width=.75pt] (g0) -- (g2);

  \node[draw, circle, scale=0.4, fill=black] (h0) at (6.5,0) {};

  \node[draw, circle, scale=0.4, fill=black] (h1) at (5.8,1) {};
    \node[draw, circle, scale=0.4, fill=black] (h2) at (7.2,1) {};
 
     \node[draw, circle, scale=0.4, fill=black] (h3) at (5.8,2) {};
  
 \node at (6.5,-0.5) {$\mathfrak{H}$}; 
         
\draw[line width=.75pt] (h0) -- (h1);
\draw[line width=.75pt] (h0) -- (h2);
\draw[line width=.75pt] (h1) -- (h3);

  \end{tikzpicture}
  \caption{}
\label{frames_fgh}
\end{center}
\end{figure}

Recall also that  the classical logic $\CPC = \IPC \oplus \chi(\mathfrak{F})$, the Jankov's logic $\mathsf{KC} =  \IPC \oplus \chi(\mathfrak{G})$ and the Scott logic $\mathsf{SL} = \IPC \oplus \chi(\mathfrak{H})$,
where $\mathfrak{F}$, $\mathfrak{G}$ and $\mathfrak{H}$ are the frames drawn in Figure \ref{frames_fgh}, see, e.g.,  \cite[Chapter 9]{ZhaCha_ml}.
Then 
\[\CPC = \IPC \oplus \chi(\mathfrak{F}) = \IPC \oplus \psi_{w_2} = \IPC \oplus (\neg\neg p \to p).\] 
Jankov's logic \[\mathsf{KC} =  \IPC \oplus \chi(\mathfrak{G}) = \IPC \oplus \psi_{w_3} = \IPC \oplus  (\neg \neg p\to p) \to (p\vee \neg p).\] 
However, note also that as $\mathfrak{G}$ is a p-morphic image of $\mathfrak{H}$ we have 
 $ \IPC \oplus \chi(\mathfrak{G}) =$
 \[\IPC \oplus \chi(\mathfrak{G}) \oplus \chi(\mathfrak{H})  = \IPC \oplus \psi_{w_3} \oplus \psi_{w_4} =   \IPC \oplus \psi_{w_3} \wedge \psi_{w_4} = \IPC \oplus  \neg p\vee\neg \neg p.\] 
Scott's logic $\mathsf{SL} =$ 
\[\IPC \oplus \chi(\mathfrak{H}) = \IPC \oplus \psi_{w_4} = \IPC \oplus g_5(p) = \IPC \oplus  ((\neg\neg p\to p)\to p\vee\neg p)\to \neg p\vee\neg\neg p.\]
 However, not all intermediate logics are axiomatizable by de Jongh formulas (see e.g., \cite[Section 9.5]{ZhaCha_ml}). We refer to  \cite[Chapter 9]{ZhaCha_ml} and \cite{BezhJank} for a detailed account of generalizations of Jankov-de Jongh formulas that axiomatize all intermediate logics. 


 Note that given a frame $\mathfrak{F}$, the number of variables used in the de Jongh formula of $\mathfrak{F}$ is the number $n$ such that there is a valuation $V$ on $\mathfrak{F}$ with
$(\mathfrak{F}, V)$ being a generated submodel of $\mathcal{U}(n)$. In the Jankov formula of $\mathfrak{F}$  the number of variables coincides with the number of elements of the Heyting algebra of all upsets of $\mathfrak{F}$;
and in the Fine-Zakharyaschev approach the number of variables of $\chi(\mathfrak{F})$ equals the cardinality of $\mathfrak{F}$. Thus, in general, de Jongh formulas need fewer varibales than Jankov and Fine-Zakharyaschev formulas. 
Often the number of variables used in such formulas is not that relevant. 
After all we are interested in semantic properties of these  formulas and what is important is that one can write such formulas and it is less important exactly how. 
However,  the syntactic shape of these and related formulas also  plays a role. We refer to \cite{NickDick18} and \cite{DickFanJulia21} for connections of subframe and stable formulas, defined in de Jongh-style, with $\mathsf{NNIL}$ 
and $\mathsf{ONNILLI}$ formulas and 
to \cite{Lauridsen19, LauridsenPhD} for their connection with the low levels of the so--called substructural  hierarchy.
Below we demonstrate that the number of variables used in the Jankov-de Jongh formula is also essential for dealing with some particular problems.

We address the following question: Given a logic $L$,  when does a logic $L ' = L \oplus  \chi(\mathfrak{F})$ have the finite model property (FMP)? In general, it is still an open problem when $\mathsf{IPC} \oplus \chi(\mathfrak{F})$ or 
in the modal case $\mathsf{S4} \oplus \chi(\mathfrak{F})$ have the FMP \cite{Kracht93}. It is known however that $\mathsf{IPC} \oplus \alpha$ and $\mathsf{S4} \oplus \alpha$ have the FMP if $\alpha$ is a formula in one variable \cite[Theorem 11.58 and Corollary 11.59]{ZhaCha_ml}.
Using the Jankov or Fine-Zakharyaschev method we obtain that the $\chi(\mathfrak{F})$ is in one variable only if $|\mathfrak{F}|\leq 2$. Hence this result cannot be applied to these formulas. 
However, using de Jongh's method (and its modal analogue) one can show that if a rooted $\mathfrak{F}$ is such that there is $V$ with $(\mathfrak{F},V)$ being a generated submodel of $\mathcal{U}(1)$
or respectively of $\mathcal{U}_{S4}(1)$ (the $1$-universal model of $\mathsf{S4}$) in the modal case, then  $\chi(\mathfrak{F})$ is in one variable. Consequently, $\mathsf{IPC} \oplus \chi(\mathfrak{F})$ and $\mathsf{S4} \oplus \chi(\mathfrak{F})$ have the FMP. 
This is more interesting in the case of  $\mathsf{S4}$, since many models of width 2 can be embedded into $\mathcal{U}_{S4}(1)$.
For example, the two frames $\mathfrak{F}_1$ and $\mathfrak{F}_2$ drawn in Figure \ref{frames_s4} admit valuations such that the corresponding models are  (isomorphic to) generated submodels of $\mathcal{U}_{S4}(1)$ \cite[Chapter 8]{ZhaCha_ml}. Therefore, 
$\mathsf{S4} \oplus \chi(\mathfrak{F}_1)$ and $\mathsf{S4} \oplus \chi(\mathfrak{F}_2)$ have the FMP.

\begin{figure}[t]
\begin{center}
\begin{tikzpicture}[scale=1, transform shape]

  \node[draw, circle, scale=0.4, fill=black] (f0) at (5,0) {};

  \node[draw, circle, scale=0.4, fill=black] (f1) at (5,1) {};
    \node[draw, circle, scale=0.4, fill=black] (f2) at (5,2) {};
      \node[draw, circle, scale=0.4, fill=black] (f3) at (5,3) {};
 \draw[line width=.75pt] (f0) -- (f3);
   
 \node at (5,-0.5) {$\mathfrak{F}_2$}; 
         
%
%
%
%
         

  \node[draw, circle, scale=0.4, fill=black] (h0) at (0,0) {};

  \node[draw, circle, scale=0.4, fill=black] (h1) at (-.7,1) {};
    \node[draw, circle, scale=0.4, fill=black] (h2) at (.7,1) {};
 
     \node[draw, circle, scale=0.4, fill=black] (h3) at (-.7,2) {};
         \node[draw, circle, scale=0.4, fill=black] (h4) at (.7,2) {};
     \node[draw, circle, scale=0.4, fill=black] (h5) at (-.7,3) {};
         \node[draw, circle, scale=0.4, fill=black] (h6) at (.7,3) {};
     \node[draw, circle, scale=0.4, fill=black] (h7) at (-.7,4) {};
         \node[draw, circle, scale=0.4, fill=black] (h8) at (.7,4) {};  
 \node at (0,-0.5) {$\mathfrak{F}_1$}; 
         
\draw[line width=.75pt] (h0) -- (h1);
\draw[line width=.75pt] (h0) -- (h2);
\draw[line width=.75pt] (-.7,1) -- (-.7,4);
\draw[line width=.75pt] (.7,1) -- (.7,4);
\draw[line width=.75pt] (h5) -- (h4);
\draw[line width=.75pt] (h6) -- (h3);
\draw[line width=.75pt] (h3) -- (h2);
\draw[line width=.75pt] (h4) -- (h1);
\draw[line width=.75pt] (h8) -- (h5);
\draw[line width=.75pt] (h7) -- (h6);

  \end{tikzpicture}
  \caption{}
\label{frames_s4}
\end{center}
\end{figure}

%




\section{Team-based intuitionistic logic and de Jongh formulas
}
\label{sec:standard_team}


In the rest of the paper, we  investigate intuitionistic and intermediate logics as well as de Jongh formulas in the setting of team semantics.  Our starting point is the intuitionistic logic over team semantics (\TIPC) introduced in \cite{CiardelliIemhoffYang20}. We show in this section, by using a disjunctive normal form of \TIPC, that universal models $\mathcal{U}(n)$ for intuitionistic logic \IPC over the standard (single-world) semantics behave also as universal models for  the intuitionistic logic \TIPC over team semantics, and thus de Jongh formulas for \TIPC can be defined in the usual manner. Furthermore,  we show that a large class of intermediate axioms (including de Jongh formulas) define the same class of Kripke frames in the team-based intuitionistic logic as in the standard (single-world-based) intuitionistic logic. This results in a first approach to obtaining intermediate logics in the team semantics setting. We will explore an alternative approach in the next section.

Let us start by recalling the definition and basic properties of team-based intuitionistic propositional logic as introduced in \cite{CiardelliIemhoffYang20}. The reader is referred to \cite{CiardelliIemhoffYang20} for further discussions. We extend the syntax of intuitionistic propositional logic 
 with an additional disjunction $\vvee$, called {\em global disjunction}; the other disjunction $\vee$ is thus referred to as {\em local disjunction}. We shall call in this paper the resulting language, {denoted as $[\bot,\wedge,\vee,\vvee,\to]$}, {\em team language} or {\em intuitionistic team logic}.  
 To be precise,  formulas of  team language are formed by the following grammar:
\[\phi::=p\mid\bot\mid\phi\wedge\phi\mid\phi\vee\phi\mid\phi\vvee\phi\mid\phi\to\phi.\]
Formulas in the standard language of intuitionistic logic (i.e., $\vvee$-free  formulas), {denoted as  $[\bot,\wedge,\vee,\to]$}, are referred to as {\em standard formulas}. 
In the rest of the paper, we reserve the first Greek letters $\alpha,\beta,\gamma,\dots$ for standard formulas, and the last Greek letters $\phi,\psi,\chi,\dots$ stand for arbitrary formulas in the team language.

Formulas of intuitionistic team logic are evaluated over the usual  intuitionistic Kripke models $\mathfrak{M} = (W , R, V )$ but with respect to {\em teams}.
A {\em team} is a set $t\subseteq W$ of possible worlds, namely $t\in \wp(W)$. 
We  define a natural relation $R^\circ\subseteq \wp(W)\times\wp(W)$ between teams by lifting the relation $R$ on $W$  to the power set of $W$ as:
\begin{itemize}
\item $tR^\circ s$ ~iff~ $\forall w\in t\exists v\in s:wRv$ (i.e., $t\subseteq R^{-1}[s]$) and $\forall v\in s\exists w\in t: wRv$ (i.e., $s\subseteq R[t]$). 
\end{itemize}
We shall often abuse the notation and write simply $tRs$ for $tR^\circ s$. 

\begin{fact}\label{fact_Rcirc_prop}
Let $\mathfrak{M}=(W,R,V)$ be a model, and $t$ and $s$ teams.
\begin{enumerate}[(i)]
\item $\varnothing R\varnothing$, $tRt$, and $tR^\circ R[t]$.
\item $s\subseteq R[t]$ ~iff~ $t\supseteq \circ R s$, where $\supseteq \circ R:=\{(t,s)\mid \exists r\subseteq W: t\supseteq r\text{ and }rRs\}$ denotes the composition of $\supseteq$ and $R$.
\end{enumerate}
\end{fact}
\begin{proof}
Item (i) is immediate. For item (ii), if $s\subseteq R[t]$, let $r=R^{-1}[s]\cap t$. Clearly, $t\supseteq r Rs$. Conversely, if $t\supseteq r R s$, then $s\subseteq R[r]\subseteq R[t]$.
\end{proof}

\begin{fact}\label{fact_Rcirc_po}
Let $\mathfrak{M}=(W,R,V)$ be a model. We have that
$\supseteq \circ R^\circ$ is a pre-order.
\end{fact}
\begin{proof}
Reflexivity follows from the fact that $t\supseteq tRt$. For transitivity, if $t\supseteq\circ Rs$ and $s\supseteq\circ Rr$, it is easy to verify that $t\supseteq (R^{-1}[r]\cap t) R^\circ r$. 
\end{proof}

Now, we give the definition of the team semantics for formulas in the team language, where our semantic clauses for implication $\to$ and local disjunction $\vee$ deviate slightly from the original definition presented in \cite{CiardelliIemhoffYang20}. The two alternative versions of the semantics are easily seen to be equivalent (by Fact \ref{fact_Rcirc_prop}(ii) above and Corollary \ref{local_v_alternative_df} below).

\begin{definition}\label{df_ts_standard}
Let $\mathfrak{M} = ( W , R, V)$ be an intuitionistic Kripke model and $t\subseteq W$ a team. 
We define the satisfaction relation $\mathfrak{M},t\models_{\TIPL}\phi$ (or simply $\mathfrak{M},t\models\phi$) inductively as follows:
\begin{itemize}
\item $\mathfrak{M},t\models p$ ~~iff~~ $t\subseteq V(p)$
\item $\mathfrak{M},t\models \bot$ ~~iff~~ $t=\varnothing$
\item $\mathfrak{M},t\models \phi\wedge\psi$ ~~iff~~ $\mathfrak{M},t\models \phi$ and $\mathfrak{M},t\models \psi$
\item $\mathfrak{M},t\models \phi\vee\psi$ ~~iff~~ there are $s,r\subseteq W$ such that $t\subseteq s\cup r$, $\mathfrak{M},s\models \phi$ and $\mathfrak{M},r\models \psi$
\item $\mathfrak{M},t\models \phi\vvee\psi$ ~~iff~~ $\mathfrak{M},t\models \phi$ or $\mathfrak{M},t\models \psi$
\item $\mathfrak{M},t\models \phi\to\psi$ ~~iff~~ for all $s\subseteq W$ with $t\supseteq \circ R s$, $\mathfrak{M},s\models \phi$ implies $\mathfrak{M},s\models \psi$
\end{itemize}
\end{definition}
We write $\mathfrak{M}\models_{\TIPL}\phi$ if $\mathfrak{M},t\models\phi$ for all teams $t\subseteq W$, and write $\mathfrak{F}\models_{\TIPL}\phi$ if $(\mathfrak{F},V)\models_{\TIPL}\phi$ for all valuations $V$ on $\mathfrak{F}$. As usual, for any class $\mathsf{F}$ of frames, we write $\mathsf{F}\models_{\TIPL}\phi$ if $\mathfrak{F}\models_{\TIPL}\phi$ for all $\mathfrak{F}\in\mathsf{F}$.
We write $\Gamma\models_{\TIPL}\phi$ (or simply $\Gamma\models\phi$) if for all Kripke models $\mathfrak{M}$ and all teams $t$, $\mathfrak{M},t\models \Gamma$ implies $\mathfrak{M},t\models \phi$.


A formula $\phi$ is said to be {\em flat}, if for all models $\MM$ and all teams $t$,
\[\mathfrak{M},t\models\phi\iff \forall w\in t: \mathfrak{M},\{w\}\models\phi\tag{\textsf{Flatness Property}}\]
For standard formulas the above defined team semantics essentially reduces to the standard (single-world) Kripke semantics, as standard formulas are flat. 

\begin{lemma}[Flatness of standard formulas]\label{flat_ipl}
Standard formulas  are flat.
As a consequence, for any set $\Delta\cup\{\alpha\}$ of  standard formulas, 
\[\Delta\models_{\TIPL}\alpha \iff \Delta\vdash_{\IPC}\alpha.\]
\end{lemma}
\begin{proof}
The flatness of standard formula is proved by a routine verification by induction. The second part of the lemma follows from completeness of \IPC and the observation that
\begin{equation}\label{flat_ipl_eq1}
\mathfrak{M},\{w\}\models_{\TIPL}\alpha\iff \mathfrak{M},w\models\alpha
\end{equation}
holds for all standard formulas $\alpha$.
\end{proof}


As with the standard intuitionistic logic,  formulas of intuitionistic team logic are also persistent, but in a more general sense. 

\begin{lemma}[Persistence]\label{persistence_TIPC} Let $\mathfrak{M}$ be a model, $t$ and $s$ teams, and $\phi$ a formula.
 If $\mathfrak{M},t\models\phi$ and $t\supseteq \circ R s$, 
 then $\mathfrak{M},s\models\phi$. In particular, every formula $\phi$ is {\em downward closed}, that is, 
 \[\mathfrak{M},t\models\phi\text{ and }s\subseteq t\Longrightarrow \mathfrak{M},s\models\phi\tag{\textsf{Downward Closure Property}}\]
\end{lemma}
\begin{proof}
By a routine inductive argument.
\end{proof}

As an immediate consequence of the downward closure property, the semantic clause of the local disjunction $\vee$ can be alternatively formulated as follows (which is the original formulation of the clause in  \cite{CiardelliIemhoffYang20}).

\begin{corollary}\label{local_v_alternative_df}
For any model $\mathfrak{M}$, team $t$ and any formulas $\phi$ and $\psi$, we have that
\[\mathfrak{M},t\models \phi\vee\psi\iff \exists s,r\subseteq t\text{ s.t. }t= s\cup r, ~\mathfrak{M},s\models \phi\text{ and }\mathfrak{M},r\models \psi.\]
\end{corollary}

The empty set $\varnothing$, being the smallest team (with respect to $\subseteq$), plays a special role in intuitionistic team logic. In particular, all formulas  have the empty team property.

\begin{lemma}[Empty team property]
For any model $\mathfrak{M}$ and any formula $\phi$, we have that
$\mathfrak{M},\varnothing\models\phi$.
\end{lemma}
\begin{proof}
A routine verification by induction.
\end{proof}

Another closely related property is the union closure property.
 A formula $\phi$ is said to be {\em closed under unions} if for any model $\MM$ and nonempty set $T$ of teams,
 \[\MM,t\models\phi\text{ for all }t\in T \Longrightarrow \MM,\bigcup T\models\phi.\tag{\textsf{Union Closure Property}}\]


\begin{fact}
A formula $\phi$ is flat iff $\phi$ satisfies the empty team, downward closure and union closure property. In particular, all standard formulas  are closed under unions.
\end{fact}
\begin{proof}
The main claim is easy to verify, and the ``in particular" part follows from \Cref{flat_ipl}. 
\end{proof}

The local disjunction $\vee$ that is inherited from the standard intuitionistic logic has the usual disjunction property in intuitionistic team logic, while the global disjunction $\vvee$ turns out to satisfy a stronger disjunction property.

\begin{lemma}[Disjunction property]\label{DP_TIPC}Let $\phi$ and $\psi$ be two arbitrary formulas.
\begin{enumerate}[(i)]
\item If $\models\phi\vee\psi$, then $\models\phi$ or $\models\psi$.
\item For any set $\Delta$ of  standard formulas, 
\[\Delta\models\phi\vvee\psi\Longrightarrow \Delta\models\phi\text{ or }\Delta\models\psi.\]
\end{enumerate}
\end{lemma}
\begin{proof}
We refer to \cite{CiardelliIemhoffYang20} for a detailed proof. Here we only sketch a proof for item (ii). Suppose $\Delta\not\models\phi$ and $\Delta\not\models\psi$. Then there are models $\mathfrak{M}_1,\mathfrak{M}_2$ and teams $t_1,t_2$ such that 
\[\mathfrak{M}_1,t_1\models\Delta,~~\mathfrak{M}_1,t_1\not\models\phi,~~\mathfrak{M}_2,t_2\models\Delta,~~\text{ and }~~\mathfrak{M}_2,t_2\not\models\chi.\]
Consider the disjoint union $\mathfrak{M}=\mathfrak{M}_1\uplus\mathfrak{M}_2$ of $\MM_1$ and $\MM_2$. It is easy to show that  $\mathfrak{M},t_1\models\Delta$, $\mathfrak{M},t_1\not\models\phi$, $\mathfrak{M},t_2\models\Delta$ and $\mathfrak{M},t_2\not\models\chi$. Since formulas in $\Delta$ are flat (by Lemma \ref{flat_ipl}), we have also that $\mathfrak{M},t_1\cup t_2\models\Delta$. Moreover, by persistence, we also obtain $\mathfrak{M},t_1\cup t_2\not\models\phi$ and $\mathfrak{M},t_1\cup t_2\not\models\psi$. Thus, $\mathfrak{M},t_1\cup t_2\not\models\phi\vvee\psi$. Hence $\Delta\not\models\phi\vvee\psi$.
\end{proof}

Closely related to the disjunction property is the fact that for any standard formula $\alpha$, the following formula (known as the {\em \textsf{Split} axiom}) is a validity:
\[(\alpha\to \phi\vvee\psi)\to(\alpha\to\phi)\vvee(\alpha\to\psi).\tag{\textsf{Split}}\]
The \textsf{Split} axiom is often considered the distinguishing axiom of team-based propositional logics. In particular, by using this validity, it is easy to show that every formula in the team language can be transformed into an equivalent formula in disjunctive normal form.

\begin{theorem}[Disjunctive normal form]\label{DNF_TIPL}
For any formula $\phi$, we have that
\(\phi\equiv\bigvvee_{i\in I}\alpha_i\)
for some (finite) set  $\{\alpha_i\mid i\in I\}$ of standard formulas.
\end{theorem}
\begin{proof}
The proof is by an easy inductive; use $\alpha\to\phi\vvee\psi\equiv(\alpha\to\phi)\vvee(\alpha\to\psi)$ for the implication case.
\end{proof}

 A sound and complete natural deduction system for intuitionistic team logic was introduced in \cite{CiardelliIemhoffYang20}, and the completeness theorem was proved by an argument that makes essential use of the above disjunctive normal form. It is easy to show that the system in \cite{CiardelliIemhoffYang20} can be presented as the following Hilbert-style system, where the formula $\alpha$ below ranges over standard formulas (i.e., $\vvee$-free formulas) only. We emphasize that the system below  is thus not closed under {\em uniform substitution}.

\begin{definition}\label{tipc_df}
The Hilbert system for intuitionistic team logic (denoted as \TIPC) consists of the following axioms:
\begin{enumerate}[(1)]
\item $\phi\to\phi\vee\psi$
\item $(\phi\to\alpha)\to ((\psi\to\alpha)\to(\phi\vee\psi\to\alpha))$  
\item $(\phi\to\chi)\to(\phi\vee\psi\to\chi\vee\psi)$
\item $\phi\vee\psi\to\psi\vee\phi$
\item $(\phi\vee\psi)\vee\chi\to\phi\vee(\psi\vee\chi)$
\item $\phi\vee(\psi\vvee\chi)\to (\phi\vee\psi)\vvee(\phi\vee\chi)$ 
\item All $\IPC$ axioms for $\bot,\wedge,\vvee,\to$, namely,
 \begin{itemize}
\item $\phi\to(\psi\to\phi)$
\item $\big(\phi\to(\psi\to\chi)\big)\to\big((\phi\to\psi)\to (\phi\to\chi)\big)$
\item $(\phi\wedge\psi)\to\phi$, $(\phi\wedge\psi)\to\psi$
\item $\phi\to(\psi\to(\phi\wedge\psi))$
\item $\phi\to(\phi\vvee\psi)$, $\psi\to(\phi\vvee\psi)$
\item $(\phi\to\chi)\to\big((\psi\to\chi)\to((\phi\vvee\psi)\to\chi)\big)$
\item $\bot\to\phi$
\end{itemize}

\item $(\alpha\to\phi\vvee\psi)\to(\alpha\to\phi)\vvee(\alpha\to\psi)$ 
\end{enumerate}
and the Modus Ponens rule.
\end{definition}

Note that the axioms (4)-(6) and (8) can also be written using bi-implications ($\leftrightarrow$) as the main connectives instead of implications ($\to$), since the right to left direction of these implications follows from the other axioms in the system (e.g., $(\phi\vee\psi)\vvee(\phi\vee\chi)\to\phi\vee(\psi\vvee\chi)$ follows from axiom (3)).


When restricted to standard formulas, the above system \TIPC coincides with \IPC; in particular, the  system \TIPC is closed under {\em standard substitutions}, i.e., substitutions that map propositional variables to standard formulas.  
The logic \TIPC is thus a conservative extension of \IPC in the following sense.



\begin{corollary}\label{IPL_flat}
For any set $\Delta\cup\{\alpha\}$ of standard formulas,
\[
\Delta\vdash_{\TIPC}\alpha\iff\Delta\vdash_{\IPC}\alpha.
\]
\end{corollary}
\begin{proof}
The right to left direction is clear; the other direction  follows easily from 
Lemma \ref{flat_ipl}.
\end{proof}

As an interesting application of the disjunctive normal form of intuitionistic team logic, we show next that the  universal models $\mathcal{U}(n)$ of \IPC behave  as universal models also for \TIPC.


\begin{proposition}
For any $n\in\mathbb{N}$, let $\mathcal{U}(n)$ be the $n$-universal model of \IPC. For any formula $\phi$, we have that
\[\mathcal{U}(n)\models_{\TIPL}\phi\iff \vdash_{\TIPC}\phi.\]
\end{proposition}
\begin{proof}
It suffices to prove the direction ``$\Longrightarrow$". Suppose $\mathcal{U}(n)\models_{\TIPL}\phi$. By Theorem \ref{DNF_TIPL}, we have that $\phi\equiv \alpha_1\vvee\dots\vvee\alpha_k$ for some standard formulas $\alpha_1,\dots,\alpha_k$. Thus, $\mathcal{U}(n)\models_{\TIPL}\alpha_i$ for some $1\leq i\leq k$, which means in particular that $\mathcal{U}(n),W\models_{\TIPL}\alpha_i$ for the domain $W$ of $\mathcal{U}(n)$. Since $\alpha_i$ is a standard formula, by flatness, we have that  $\mathcal{U}(n),w\models\alpha_i$ for all $w\in W$. Now, by Theorem \ref{Un_prop}, we obtain that $\vdash_{\IPC}\alpha_i$, which implies by Corollary \ref{IPL_flat} that $\vdash_{\TIPC}\alpha_i$. Hence $\vdash_{\TIPC}\alpha_1\vvee\dots\vvee\alpha_k$, thereby $\vdash_{\TIPC}\phi$.
\end{proof}

Having the universal models in place for \TIPC, we can consider the  de Jongh formulas also in the context of \TIPC. Thanks again to the disjunctive normal form of \TIPC, a large class of the same de Jongh formulas and intermediate axioms for the standard language of \IPC still characterize the same class of frames in the setting of \TIPC. Before we present our result, let us first define intermediate logics 
in this setting. This definition has appeared also in a recent work \cite{Quadrellaro21}; cf.~also a similar but subtly different definition in \cite{Puncochar21}.

\begin{definition}\label{team_intermediate_df1}
Let $\IPC\subseteq \LL\subseteq \CPC$ be an intermediate logic with $\LL=\IPC\oplus\Delta$ and $\Delta=\{\alpha_i\mid i\in I\}$ a set of axioms written in the standard language (i.e., each $
\alpha_i$ is a formula in  $[\bot,\wedge,\vee,\to]$). Denote by $\mathsf{t}\LL$ the logic obtained by closing the set
\[\TIPC\cup\{\alpha_i(\vec{\beta}/\vec{p})\mid i\in I,~\text{and each }\beta_j\text{ is a standard formula}\}\]  
under Modus Ponens. 
We write $\vdash_{\mathsf{t}\LL}\phi$ if $\phi\in \mathsf{t}\LL$. Note that $\vdash_{\mathsf{t}\LL}\phi$ iff $\Delta\vdash_{\TIPC}\phi$.
\end{definition}

For instance, consider Jankov's logic  $\mathsf{KC}=\IPC \oplus \chi(\mathfrak{G})=\IPC\oplus \neg p\vee\neg\neg p$, where the frame $\mathfrak{G}$ is as drawn in Figure \ref{frames_fgh}. 
The intermediate logic $\tl\mathsf{KC}$ is then obtained by closing the set
\[\TIPC\cup\{\neg\alpha\vee\neg\neg \alpha\mid \alpha\in [\bot,\wedge,\vee,\to]\}\]
under Modus Ponens.  Note that $\tl\mathsf{KC}$ may not be  closed under uniform substitution, 
as, e.g., an arbitrary formula $\neg\phi\vee\neg\neg\phi$ (with global disjunction $\vvee$) is not necessarily  in $\mathsf{tKC}$.

The logic $\tl\mathsf{CPC}$ obtained by extending $\TIPC$ with the classical axiom $\neg\neg p\to p$ 
is known in the literature as a variant of {\em propositional dependence logic} (see \cite{Ciardelli2015,VY_PD}), and its fragment with connectives $\bot,\wedge,\vvee,\to$, {denoted as $\mathop{[\bot,\wedge,\vvee,\to]}$}, is called {\em inquisitive logic} (\cite{InquiLog}). We emphasize again that the logic $\tl\mathsf{CPC}$ is not closed under uniform substitution; in particular its $[\bot,\wedge,\vvee,\to]$-fragment, i.e., inquisitive logic (which contains all \IPC axioms and the axiom $\neg\neg\alpha\to\alpha$) is not equivalent to classical logic, as an arbitrary formula $\neg\neg\phi\to\phi$ is not necessarily in the logic.
It is interesting to note that for inquisitive logic, or the $[\bot,\wedge,\vvee,\to]$-fragment of $\tl\mathsf{CPC}$, the \textsf{Split} axiom is equivalent to the Kreisel-Putnam  axiom
\[(\neg p\to q\vvee r)\to (\neg p\to q)\vvee(\neg p\to r)\tag{\textsf{KP}}.\]
The logic can thus be viewed as the intermediate theory (in the language $\mathop{[\bot,\wedge,\vvee,\to]}$)  
obtained by extending the Kreisel-Putnam logic $\textsf{KP}$ (i.e., $\IPC\oplus\textsf{KP}$) with the double negation elimination axiom $\neg\neg p\to p$. It was further shown that the logic is equivalent to the so-called {\em negative variant} of Medvedev's logic \textsf{ML} of finite problems (whose frames are of the form $(\wp(W)\setminus\{\varnothing\},\supseteq)$ for finite sets $W$). The reader is referred to \cite{InquiLog} for detailed discussions. In Section \ref{sec:pw} we will explore a related 
but different power set construction than that of the Medvedev frames. 



Now we are ready to show the main result of this section that standard axioms $\Delta$ 
define the same classes of frames in \TIPC as in \IPC, as long as the intermediate logic $\IPC\oplus\Delta$  has the disjunction property or is canonical.

\begin{theorem}\label{inter_compl_1}
Let $\LL=\IPC\oplus\Delta$ be an intermediate logic that  is complete with respect to a class $\mathsf{F}_{\LL}$ of $\LL$-frames, where $\Delta$ is a set of standard formulas. If $\LL$ has the disjunction property or is canonical, then for every formula $\phi$, we have that
\[\vdash_{\tl\LL}\phi\iff \mathsf{F}_{\mathsf{L}}\models_{\TIPL}\phi.\]
%
\end{theorem}
\begin{proof}
By Theorem \ref{DNF_TIPL} and the completeness theorem of \TIPC, we may assume that $\phi=\alpha_1\vvee\dots\vvee\alpha_k$ for some standard formulas $\alpha_1,\dots,\alpha_k$. Suppose that $\vdash_{\tl\LL}\phi$, i.e., $\Delta\vdash_{\TIPC}\alpha_1\vvee\dots\vvee\alpha_k$. By the disjunction property of $\vvee$ (by \Cref{DP_TIPC} and completeness of \TIPC), we have that $\Delta\vdash_{\TIPC}\alpha_i$  for some $1\leq i\leq k$. It then follows from Corollary \ref{IPL_flat} that $\Delta\vdash_{\IPC}\alpha_i$, which implies $\mathsf{F}_{\mathsf{L}}\models\alpha_i$ in the sense of the standard (single world-based) semantics of \IPC. By flatness of $\alpha_i$ (\Cref{flat_ipl}), we obtain $ \mathsf{F}_{\mathsf{L}}\models_{\TIPL}\alpha_i$ as well. Hence, $\mathsf{F}_{\mathsf{L}}\models_{\TIPL}\alpha_1\vvee\dots\vvee\alpha_k$, namely $\mathsf{F}_{\mathsf{L}}\models_{\TIPL}\phi$.

For the converse direction, we first treat the case when $\mathsf{L}$ has the disjunction property. Suppose $\mathsf{F}_{\mathsf{L}}\models_{\TIPL}\alpha_1\vvee\dots\vvee\alpha_k$. Since $\phi\vvee\psi\models_{\TIPL}\phi\vee\psi$ (by empty team property), we have $\mathsf{F}_{\mathsf{L}}\models_{\TIPL}\alpha_1\vee\dots\vee\alpha_k$ as well, where the standard formula $\alpha_1\vee\dots\vee\alpha_k$ is flat (by \Cref{flat_ipl}). Thus, $\mathsf{F}_{\mathsf{L}}\models\alpha_1\vee\dots\vee\alpha_k$ in the sense of the standard (single world-based) semantics of \IPC, which implies further that $\Delta\vdash_{\IPC}\alpha_1\vee\dots\vee\alpha_k$. Since $\mathsf{L}$ has the disjunction property, we have $\Delta\vdash_{\IPC}\alpha_i$ for some $1\leq i\leq k$. Hence, by Corollary \ref{IPL_flat}, we conclude that $\Delta\vdash_{\TIPC}\alpha_i$, which gives $\Delta\vdash_{\TIPC}\alpha_1\vvee\dots\vvee\alpha_k$, i.e., $\vdash_{\tl\LL}\phi$.

Next, we treat the case when $\mathsf{L}$ is canonical. Suppose that $\Delta\nvdash_{\TIPC}\alpha_1\vvee\dots\vvee\alpha_k$. Then, for all $1\leq i\leq k$, we have that $\Delta\nvdash_{\TIPC}\alpha_i$, which implies $\Delta\not\vdash_{\IPC}\alpha_i$ by Corollary \ref{IPL_flat}. Now, for each $1\leq i\leq k$,  the canonical model $ \mathfrak{M}_{\mathsf{L}}^c$ of $\mathsf{L}$ contains a witness $w_i$ satisfying $\mathfrak{M}_{\mathsf{L}}^c,w_i\not\models\alpha_i$. Thus, $\mathfrak{M}_{\mathsf{L}}^c,\{w_i\}\not\models_{\TIPL}\alpha_i$ by Expression (\ref{flat_ipl_eq1}). Let $t=\{w_1,\dots,w_k\}$. By persistence (\Cref{persistence_TIPC}), we must have that $\mathfrak{M}_{\mathsf{L}}^c,t\not\models_{\TIPL}\alpha_i$ for all $1\leq i\leq k$, and thus $\mathfrak{M}_{\mathsf{L}}^c,t\not\models_{\TIPL}\alpha_1\vvee\dots\vvee\alpha_k$. Since $\mathsf{L}$ is canonical, we know that $\mathfrak{M}_{\mathsf{L}}^c\in\mathsf{F}_{\mathsf{L}}$, and thus $\mathsf{F}_{\mathsf{L}}\not\models_{\TIPL}\alpha_1\vvee\dots\vvee\alpha_k$.
\end{proof}

Let us give some examples of the applications of the above theorem. Since Jankov's logic $\mathsf{KC} =  \IPC \oplus \chi(\mathfrak{G})$  is  canonical, by the above theorem, 
\[
\vdash_{\tl\mathsf{KC}}\phi\iff\mathsf{F}_{\mathsf{KC}}\models_{\TIPL}\phi,\]
where $\mathsf{F}_{\mathsf{KC}}$ is the class of finite rooted frames with unique top points. 
Since  classical propositional logic $\CPC = \IPC \oplus \chi(\mathfrak{F})$ (where the frame $\mathfrak{F}$ is as drawn in  Figure \ref{frames_fgh}) is canonical, the logic $\tl\mathsf{CPC}$ 
is complete with respect to all classical frames (i.e., frames with dead points only). Note that over classical Kripke models, a team (i.e., a set of possible worlds with no proper successors) can be identified as a set of propositional valuations. In fact, the teams for the logic $\tl\mathsf{CPC}$ (or propositional dependence logic \cite{VY_PD} and inquisitive logic \cite{InquiLog}) were originally introduced as such sets of valuations.

Theorem \ref{inter_compl_1}, however, does not apply for the intermediate logics that neither have the disjunction property nor are canonical. For instance, Scott's logic $\mathsf{SL} =\IPC \oplus \chi(\mathfrak{H})$ (where the frame $\mathfrak{H}$ is as drawn in  Figure \ref{frames_fgh}) is one of such examples. Studying the properties of the corresponding intermediate team logics for such logics is left as future work.

\section{Generalized team intuitionistic Kripke semantics}\label{sec:gt}

In the previous section we have defined (in Definition \ref{team_intermediate_df1}) intermediate team logics by adding to \TIPC a set of standard axioms with the local disjunction $\vee$. As shown in Theorem \ref{inter_compl_1}, a large class of such axioms characterize, in the usual manner, properties of the underlying Kripke frames. One important feature of such defined intermediate team logics is that they all validate the \textsf{Split} axiom, which characterizes certain unique properties of the underlying structure of the set of teams. Roughly speaking, these logics with two disjunctions (the local one $\vee$ and the global one $\vvee$) can be considered logics with two layers: The first layer, which we shall call the {\em core} or {\em inner} layer, is characterized by the different properties of the underlying Kripke frames. By the result of Theorem \ref{inter_compl_1}, this layer can be axiomatized by standard formulas, i.e., formulas that contain the local disjunction $\vee$ only, to a large extent (i.e., for axioms that define intermediate logics with disjunction property or are canonical). The second layer, which we shall call the {\em team} or {\em outer} layer, is characterized by the properties of the underlying structure of the set of teams. In the setting of Definition \ref{team_intermediate_df1}, this team layer is kept fixed by always adopting the  \textsf{Split} axiom in the intermediate team logics (as these logics are always extensions of \TIPC). Different intermediate logics differ only in the core layer, through validating different standard axioms. See Figure \ref{fig:standard_logics} for an illustration of such constructions.

\begin{figure}[t]\caption{Intermediate team logics obtained by changing the core/team layer}\label{fig:generalized_logics}
\begin{center}
\begin{tikzpicture}[transform shape, scale=1]

	\draw[dashed,color=gray] (4,-1.1) arc (270:90:1.2 and 2.2);
	\draw[semithick] (4,-1.1) arc (-90:90:1.2 and 2.2);
\draw[semithick] (0,-1.1)  arc (270:90:1.2 and 2.2);
\draw[semithick] (0,-1.1) arc (-90:90:1.2 and 2.2);
	\draw[semithick] (4,-1.11) -- (0,-1.11);
	\draw[semithick] (4,3.3) -- (0,3.3);

	\draw[dashed,color=gray] (4,-0.6) arc (270:90:0.9 and 1.7);
	\draw[semithick] (4,-0.6) arc (-90:90:0.9 and 1.7);
\draw[semithick] (0,-0.6)  arc (270:90:0.9 and 1.7);
\draw[semithick] (0,-0.6) arc (-90:90:0.9 and 1.7);
       
	\draw[dashed,color=gray] (4,-0.2) arc (270:90:0.57 and 1.35);
	\draw[semithick] (4,-0.2) arc (-90:90:0.57 and 1.35);
\draw[semithick] (0,-0.2)  arc (270:90:0.57 and 1.35);
\draw[semithick] (0,-0.2) arc (-90:90:0.57 and 1.35);

	\draw[dashed,color=gray!50, fill=gray!20] (4,0) arc (270:90:0.4 and 1.1);
	\draw[semithick, color=gray!50] (4,0) -- (0,0.3);
	\draw[thick, color=white] (0,0.3) -- (1.1,0.22);
	\draw[semithick, color=gray!50, dashed] (0,0.3) -- (1.1,0.22);
	\draw[semithick, color=gray!50] (4,2.2) -- (0,1.7);
\draw[thick, color=white] (1.13,1.84) -- (0,1.7);
\draw[semithick, color=gray!50, dashed] (1.13,1.84) -- (0,1.7); 

	\draw[semithick, color=gray, fill=gray!20] (4,0) arc (-90:90:0.4 and 1.1);
	\draw[semithick, color=gray, fill=gray!20] (0,1) ellipse (0.3 and 0.7);
 
    \node at (0,-0.85) {\footnotesize $\vvee$};
   \node at (0,0.6) {\footnotesize $\vee$};
  \node at (0,1.1) {\footnotesize \IPC};
    \node at (0,-0.39) {\scriptsize \textsf{Split}};

 \node at (4,1.1) {\footnotesize \CPC};
    \node at (4,-0.39) {\scriptsize \textsf{Split}};
   \node at (4,0.6) {\footnotesize $\vee$};
      \node at (4,-0.85) {\footnotesize $\vvee$};

\end{tikzpicture}
\end{center}
\end{figure}
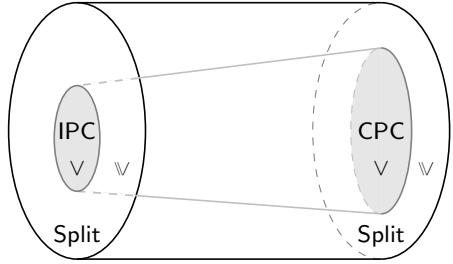

Another natural approach to define intermediate team logics is to allow also the team layer of the logic to vary, by replacing the \textsf{Split} axiom with other axioms in the full language of \TIPC (see Figure \ref{fig:generalized_logics} for an illustration of such constructions). In this section, we make a first attempt in this direction. 

Since the \textsf{Split} axiom captures the structure of the teams, defining intermediate team logics that validate other axioms  than the \textsf{Split} axiom also amounts to generalizing the original team semantics in some appropriate manner.  To this end, in Section \ref{sec:pw} we interpret the logic \TIPC in the so-called (full) powerset model and compare the logic with a variant of intuitionistic modal logic \IK. On the basis of this connection, in Section \ref{sec:gtipl}, we define generalized team intuitionistic Kripke semantics and discuss the failure of the \textsf{Split} axiom and two other axioms involving the local disjunction $\vee$ in this new setting. 
In Section \ref{sec:recover}, we provide frame conditions validating these three axioms and thus recovering \TIPC.


\begin{figure}[t]\caption{Intermediate team logics obtained by changing the core}\label{fig:standard_logics}
\begin{center}
\begin{tikzpicture}[transform shape, scale=1]


	\draw[dashed,color=gray] (4,-0.6) arc (270:90:0.9 and 1.7);
	\draw[semithick] (4,-0.6) arc (-90:90:0.9 and 1.7);
\draw[semithick] (0,-0.6)  arc (270:90:0.9 and 1.7);
\draw[semithick] (0,-0.6) arc (-90:90:0.9 and 1.7);
	\draw[semithick] (4,-0.6) -- (0,-0.6);
	\draw[semithick] (4,2.8) -- (0,2.8);

	\draw[dashed,color=gray!50, fill=gray!20] (4,0) arc (270:90:0.4 and 1.1);
	\draw[semithick, color=gray!50] (4,0) -- (0,0.3);
	\draw[thick, color=white] (0,0.3) -- (0.76,0.24);
	\draw[semithick, color=gray!50, dashed] (0,0.3) -- (0.76,0.24);
	\draw[semithick, color=gray!50] (4,2.2) -- (0,1.7);
\draw[thick, color=white] (.81,1.8) -- (0,1.7);
\draw[semithick, color=gray!50, dashed] (.81,1.8) -- (0,1.7); 

	\draw[semithick, color=gray, fill=gray!20] (4,0) arc (-90:90:0.4 and 1.1);
	\draw[semithick, color=gray, fill=gray!20] (0,1) ellipse (0.3 and 0.7);
 
    \node at (0.6,0.6) {\footnotesize $\vvee$};
   \node at (0,0.6) {\footnotesize $\vee$};
  \node at (0,1.1) {\footnotesize \IPC};
    \node at (0,-0.3) {\footnotesize \textsf{Split}};
    
      \node at (4,1.1) {\footnotesize \CPC};
    \node at (4,-0.3) {\footnotesize \textsf{Split}};
   \node at (4,0.6) {\footnotesize $\vee$};
      \node at (4.6,0.6) {\footnotesize $\vvee$};

\end{tikzpicture}
\end{center}
\end{figure}

\subsection{The power set model and intuitionistic modal logic}\label{sec:pw}


In this section, we relate  \TIPC to {\em inuitionistic modal logic} \IK (studied in, e.g., \cite{IML_PlotkinStirling86,WolterZakhIML99}) 
and give a single-world semantics interpretation of the team-based logic \TIPC. The key observation for this connection is that the team-based logic \TIPC can be given a single-world semantics over the powerset models induced by the underlying intuitionistic Kripke models. We show that the local disjunction $\vee$-free fragment of  \TIPC can be translated, in a truth-preserving manner, into the diamond-free fragment of
intuitionistic modal logic \IK over powerset models, while full \TIPC corresponds to a variant of \IK.  Our results are inspired by a similar connection and powerset models discussed in \cite{Yang17MD} in the context of modal dependence logics. The logic \TIPC can be translated into $\mathsf{S4}$ modal dependence logic via a  G\"{o}del-style translation, as shown in  \cite{CiardelliIemhoffYang20}.

We first give the translation of the local disjunction $\vee$-free fragment of \TIPC into the single-world semantics-based intuitionistic modal logic \IK, and then extend this result to full \TIPC.
Recall that formulas of \TIPC are evaluated in an intuitionistic Kripke model $\mathfrak{M}=(W,R,V)$ over {\em teams} (or sets of worlds), which are elements in the power set of $W$. Every intuitionistic Kripke model $\mathfrak{M}$ is thus naturally associated with a powerset model defined as follows. Similar powerset models were considered also in \cite{vanBenthem21}. 

  \begin{definition}\label{pwset_ml_df}
Let $\mathfrak{M}=(W,R,V)$ be an intuitionistic Kripke model. The \emph{powerset  model $\MM^\circ$ induced by \MM} is a quadruple $\MM^\circ=(W^\circ,\supseteq, R^\circ,V^\circ)$, where
\begin{itemize}
\item $W^\circ=\wp(W)\setminus\{\varnothing\}$, i.e. $W^\circ$ consists of all nonempty teams $t\subseteq W$, 
\item $\supseteq$ is the superset relation, 
\item $R^\circ$ is the lifting of $R$ on $\wp(W)$, 
\item the valuation $V^\circ$ is defined as: \(t\in V^\circ(p) \text{~~iff~~} t\subseteq V(p).\)
\end{itemize}
\end{definition}


If $\MM$ is finite, the reduct $(W^\circ,\supseteq)=(\wp(W)\setminus\{\varnothing\},\supseteq)$ of the underlying frame of the powerset model $\mathfrak{M}^\circ$ is clearly a frame for Medvedev's logic \textsf{ML}.
The whole powerset model $\MM^\circ$ is a special case of a bi-relation intuitionistic Kripke model of intuitionistic modal logic \IK. 
We now recall the definition of such models.

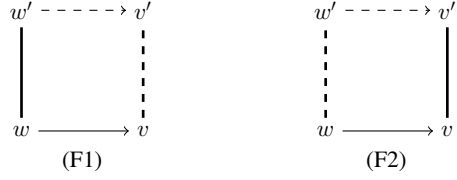
\begin{figure}[t]
 \begin{center}
\begin{tikzpicture}[scale=.8, transform shape] 
\node  (1)  at (0, 0) {\normalsize$w$};
\node  (3) at +(0,2) {\normalsize$w'$};
\node (4) at +(2,2) {\normalsize$v'$};
\node  (2) at +(2,0) {\normalsize$v$};
\draw[line width=1pt, dashed] (2) -- (4);
\draw[line width=1pt] (1) -- (3);
\draw [dashed, ->] (3) -- (4);
\draw [->] (1) -- (2);
\node at (1,-0.5) {\normalsize(F1)};
\node  (11)  at (5, 0) {\normalsize$w$};
\node  (13) at +(5,2) {\normalsize$w'$};
\node (14) at +(7,2) {\normalsize$v'$};
\node  (12) at +(7,0) {\normalsize$v$};
\draw[line width=1pt] (12) -- (14);
\draw[line width=1pt, dashed] (11) -- (13);
\draw [dashed, ->] (13) -- (14);
\draw [->] (11) -- (12);
\node at (6,-0.5) {\normalsize(F2)};


\end{tikzpicture}
 \caption{Conditions on frames. 
 The directed lines represent the $R$ relation and the undirected lines represent the $\succcurlyeq$ relation with the nodes positioned above being accessible from the ones positioned below.
 }\label{Fig_F1F2}
\end{center}
  \end{figure}

\begin{definition}\label{bi-relation_ml}
A \emph{bi-relation intuitionistic Kripke frame} 
is a triple $\mathfrak{F}=(W,\succcurlyeq,R)$, where 
\begin{itemize}
\item $W$ is a nonempty set;
\item $\succcurlyeq$ is a partial order on $W$ and $R$ is a binary relation on $W$ 
satisfying the following conditions (F1) and 
(F2) (see also Figure \ref{Fig_F1F2}): 
\begin{description}
  \item[F1:] If $w\succcurlyeq w'$ and $wR v$, then there exists $v'\in W$ such that $vR v'$ and $w'R v'$ (i.e., $\preccurlyeq\circ R\,\subseteq\, R\circ \preccurlyeq$).
  \item[F2:] If $wR v$ and $v\succcurlyeq v'$, then there exists $w'\in W$ such that $wR w'$ and $w'R v'$ (i.e., $R\circ \succcurlyeq\,\subseteq\, \succcurlyeq\circ R$).
\end{description}
\end{itemize}
A \emph{bi-relation intuitionistic Kripke model} is a pair $\mathfrak{M}=(\mathfrak{F},V)$ of a bi-relation intuitionistic Kripke frame $\mathfrak{F}=(W,\succcurlyeq,R)$
and a \emph{valuation} $V:\Prop\to\wp(W)$ satisfying \emph{persistence} with respect to $\succcurlyeq$ (i.e., $w\in V(p)$ and $w \succcurlyeq v$ imply $v\in V(p)$).
\end{definition}

It is easy to verify that the superset relation $\supseteq$ and the lifting $R^\circ$ in a powerset model $\mathfrak{M}^\circ$ satisfy conditions (F1) and (F2), and the lifting $V^\circ$ of the valuation is persistent with respect to the partial order $\supseteq$. A powerset model $\mathfrak{M}^\circ$ is thus a bi-relation intuitionistic Kripke model.
As pointed out already, such type of connection 
has already been discussed in \cite{Yang17MD}, where it was shown that modal dependence logic corresponds to a special type of \IK and the relation $R^\circ$ is the lifting of an arbitrary relation $R$ in a modal Kripke model. Here in our setting, the relation $R^\circ$ is the lifting of a partial order $R$ of an intuitionistic Kripke model.  Observe that the lifting $R^\circ$ is actually a pre-order, or a partial order  in case the model is finite.  
\begin{fact}\label{lifting_R_preorder_M0}
The lifting $R^\circ$ is a pre-order.  If $\MM$ is finite, then $R^\circ$ is a partial order.
\end{fact}
\begin{proof}
The reflexivity and transitivity of $R^\circ$ follow immediately from the reflexivity and transitivity of $R$.
If $tR^\circ s$, $sR^\circ t$ and $s,t$ are finite, then we must have $t=s$, since $R$ is anti-symmetric. 
\end{proof}

Recall that the language of intuitionistic modal logic \IK is obtained by enriching the language of intuitionistic propositional logic \IPC with modalities. That is, formulas of $\IK$ are formed by the following grammar:
\[\phi::=p\mid \bot\mid \phi\wedge\phi\mid\phi\vee\phi\mid\phi\to\phi\mid\Box\phi\mid\Diamond\phi.\]
Formulas of \IK are evaluated over bi-relation intuitionistic Kripke models with respect to single worlds.

\begin{definition}\label{sat_df_bir}
The  \emph{satisfaction} relation $\MM,w\Vdash\phi$ between a bi-relation intuitionistic Kripke model  $\MM=\mathop{(W,\succcurlyeq,R,V)}$, a world $w\in W$ and an \IK-formula $\phi$  is defined inductively as follows:
\begin{itemize}
\item $\MM,w\Vdash p$ ~~iff~~ $w\in V(p)$
\item $\MM,w\nVdash\bot$ always
\item $\MM,w\Vdash\phi\wedge\psi$ ~~iff~~ $\MM,w\Vdash\phi$ and $\MM,w\Vdash\psi$
\item $\MM,w\Vdash\phi\vee\psi$ ~~iff~~ $\MM,w\Vdash\phi$ or $\MM,w\Vdash\psi$ 
\item $\MM,w\Vdash\phi\to\psi$ ~~iff~~ for all $v\in W$ such that $w\succcurlyeq v$, if $\MM,v\Vdash\phi$, then $\MM,v\Vdash\psi$
\item $\MM,w\Vdash\Box\phi$ ~~iff~~ for all $u,v\in W$ such that $w\succcurlyeq u$ and $uR v$, it holds that $\MM,v\Vdash\phi$
\item $\MM,w\Vdash\Diamond\phi$ ~~iff~~ there exists $v\in W$ such that $wRv$ and $\MM,v\Vdash\phi$
\end{itemize}
\end{definition}

It is routine to verify that persistence of propositional variables generalizes to arbitrary \IK-formulas, namely, $\MM,w\Vdash\phi$ and $w\succcurlyeq v$ imply $\MM,v\Vdash\phi$.

Formulas in the local disjunction $\vee$-free fragment of \TIPC can be translated into the diamond-free fragment  of \IK via a truth preserving translation defined as follows. A similar translation from modal dependence logic into \IK was defined in  \cite{Yang17MD}.

\begin{definition}\label{trans_tipl2ik}
Define a translation ${(\cdot)}^\tau$ from the $\vee$-free fragment of \TIPC into (the diamond-free fragment of) \IK inductively as:
\begin{itemize}
\item $p^\tau:=p$
\item $\bot^\tau:=\bot$ 
\item $(\phi\vvee\psi)^\tau:=\phi^\tau\vee\psi^\tau$
\item $(\phi\wedge\psi)^\tau:=\phi^\tau\wedge\psi^\tau$
\item $(\phi\to\psi)^\tau:=\Box(\phi^\tau\to\psi^\tau)$
\end{itemize}
\end{definition}

We now show that the above defined translation preserves truth over powerset models and nonempty teams.

\begin{proposition}\label{comodel_model}
Let $\MM=(W,R,V)$ be an intuitionistic Kripke model and $t\subseteq W$ a nonempty team. For any $\vee$-free formula $\phi$, $\MM,t\models \phi\iff \MM^\circ,t\Vdash\phi^\tau$.
\end{proposition}

\begin{proof}
We prove the lemma by induction on $\phi$. 
If $\phi=\bot$, then for any nonempty team $t$, we have $\MM,t\not\models\bot$ and $\MM^\circ,t\nVdash\bot$.

If $\phi=p$, then $\MM,t\models p$ iff $t\subseteq V(p)$ iff $\MM^\circ,t\Vdash p$.

If $\phi=\psi\to\chi$, suppose $t\models \psi\to\chi$. We show that $\MM^\circ,t\Vdash\Box(\psi^\tau\to\chi^\tau)$. Let $r,s,u\in \wp(W)\setminus\{\varnothing\}$ be such that $t\supseteq rRs\supseteq u$ and $u\Vdash\psi^\tau$. By induction hypothesis, we have that $u\models\psi$. Since $t\supseteq rRs$ and $t\models \psi\to\chi$, by persistence (\Cref{persistence_TIPC}), we have $s\models \psi\to\chi$. Moreover, since $R$ is reflexive, we have $s\supseteq uRu$.
Thus,  we obtain $u\models \chi$, from which we conclude $u\Vdash\chi^\tau$ by induction hypothesis.

Conversely, suppose  $t\Vdash\Box(\psi^\tau\to\chi^\tau)$. Let $t\supseteq r Rs$ 
be such that $s\models\psi$. If $s=\varnothing$, then $s\models\chi$ by the empty team property. Now, assume that $s\neq\varnothing$. By induction hypothesis, we have that $s\Vdash\psi^\tau$. Moreover, since $rRs$ and $s\neq\varnothing$, we know that $r\neq\varnothing$ and $r\in W^\circ$.
 Since $t\supseteq r Rs\supseteq s$, we obtain that $s\Vdash\chi^\tau$ by assumption. Hence, $s\models\chi$ holds by induction hypothesis.




The other cases are straightforward. 
\end{proof}

Now, we turn to full \TIPC that has the local disjunction $\vee$ in the language. In the powerset model $\MM^\circ$ as defined in Definition \ref{pwset_ml_df}, we have excluded the empty team from the domain (i.e., $W^\circ=\wp(W)\setminus\{\varnothing\}$), as the empty team behaves as a trivial team in the $\vee$-free fragment of \TIPC. When it comes to full \TIPC, in the presence of $\vee$ the empty team, however, does play a special and nontrivial role. For this reason, to interpret \TIPC in the single-world semantics, we shall instead consider full powerset models $\MM^\bullet$ that has as the domain the full power set $\wp(W)$, as done also in \cite{Yang17MD} in the model dependence logic context. The local disjunction $\vee$  will be interpreted essentially as a binary diamond modality $\otimes$ (also called {\em tensor}) in the single-world semantics over the full powerset models $\MM^\bullet$. The model $\MM^\bullet$ is equipped also with a  binary operator, the set-theoretic union $\cup$, which, together with the partial order $\supseteq$, induces a ternary relation for $\otimes$. This treatment of the local disjunction is slightly different from that in \cite{Yang17MD}, where  a ternary relation for $\otimes$ was directly included in the full powerset model.


  \begin{definition}\label{fullpwset_ml_df}
Let $\mathfrak{M}=(W,R,V)$ be an intuitionistic Kripke model. The \emph{full powerset  model $\MM^\bullet$ induced by \MM} is a tuple $\MM^\bullet=(W^\bullet,\supseteq, R^\circ,\cup,\varnothing,V^\circ)$, where
\begin{itemize}
\item $W^\bullet=\wp(W)$, i.e. $W^\bullet$ consists of all teams $t\subseteq W$, 
\item $\cup$ is the set-theoretic union,
\item and the other components are as defined in the powerset model $\MM^\circ$ (Definition \ref{pwset_ml_df}).
\end{itemize}
\end{definition}



Clearly, $(W^\bullet,\cup,\varnothing)$ forms a bounded join-semilattice  (i.e., $\cup$ is  associative, commutative and idempotent, and $\varnothing\cup t=t$ for all $t\in W^\bullet$), and $\subseteq$ is the associated partial order (i.e., $\subseteq$ satisfies: $t\subseteq s$ iff $t\cup s=t$). A full powerset model can be viewed as a special case of a bi-relation intuitionistic Kripke model together with an additional binary operator $\Cup$ and a constant $0$, which we shall call a {\em tensored bi-relation intuitionistic Kripke model}.

\begin{definition}\label{tensor_bi-relation_ml}
A \emph{tensored bi-relation intuitionistic Kripke frame} 
is a quadruple $\mathfrak{F}=(W,\succcurlyeq,R,\Cup,0)$, where 
\begin{itemize}
\item $(W,\succcurlyeq,R)$ is a bi-relation intuitionistic Kripke frame;
\item $(W,\Cup,0)$ is a bounded join-semilattice; 
\item $\preccurlyeq$ is the associated partial order of $(W,\Cup,0)$. 
\end{itemize}
A \emph{tensored bi-relation intuitionistic Kripke model} is a pair $\mathfrak{M}=(\mathfrak{F},V)$ such that $\mathfrak{F}=(W,\succcurlyeq,R,\Cup,0)$ is a tensored bi-relation intuitionistic Kripke frame 
and $V$ is a \emph{valuation} satisfying \emph{persistence} with respect to $\succcurlyeq$.
\end{definition}

In order to compare the full logic \TIPC with \IK, we shall consider a (diamond-free) variant of \IK (denoted as $\IKt$) that has an additional binary diamond, the tensor $\otimes$, and a constant symbol $\mathbf{0}$. 
 We now  define the semantics for  $\IKt$ by adapting the standard semantics of \IK over bi-relation intuitionistic Kripke models. Such a language was considered already in  \cite{Yang17MD}, where the semantics for tensor $\otimes$ was defined slightly differently and the constant symbol $\mathbf{0}$ was not isolated.

\begin{definition}\label{sat_df_tensor_bir}
The  \emph{satisfaction} relation $\MM,w\Vdash\phi$ between a tensored bi-relation intuitionistic Kripke model  $\MM=\mathop{(W,\succcurlyeq,R,\Cup,0,V)}$, a world $w\in W$ and an $\IKt$-formula $\phi$   is defined inductively as in Definition \ref{sat_df_bir} except for the following cases:
\begin{itemize}
\item $\MM,w\Vdash\mathbf{0}$ ~~iff~~ $w=0$
\item $\MM,w\Vdash \phi\otimes\psi$ ~~iff~~ there are $u,v\in W$ such that $w\preccurlyeq u\Cup v$, $\MM,u\Vdash \phi$ and $\MM,v\Vdash \psi$ 
\end{itemize}
\end{definition}

We now verify that formulas of $\IKt$ are persistent over tensored bi-relation intuitionistic Kripke models.

\begin{proposition}[Persistence]\label{IK_persistence}
For any \IKt-formula $\phi$, if $\MM,w\Vdash\phi$ and $w\succcurlyeq v$, then $\MM,v\Vdash\phi$.
\end{proposition}
\begin{proof}
Let $\MM=\mathop{(W,\succcurlyeq,R,\Cup,0,V)}$. The proposition is proved by induction. If $\phi=\mathbf{0}$, then $\MM,w\Vdash\mathbf{0}$ and $w\succcurlyeq v$ imply that $w=v=0$, which then gives $\MM,v\Vdash\mathbf{0}$.

If $\phi=\psi\otimes\chi$, suppose $\MM,w\Vdash\psi\otimes\chi$ and $w\succcurlyeq v$. Then there are $x,y\in W$ such that $w\preccurlyeq x\Cup y$, $\MM,x\Vdash \phi$ and $\MM,y\Vdash \psi$. Since $v\preccurlyeq w\preccurlyeq x\Cup y$, we conclude that $\MM,v\Vdash\psi\otimes\chi$.

The other cases are standard.
\end{proof}

Formulas of \TIPC can be translated into \IKt via a truth preserving translation defined as follows.

\begin{definition}\label{trans_bold_tipl2ik}
Define  inductively a translation ${(\cdot)}^{\pmb{\tau}}$ from \TIPC into \IKt as:
\begin{itemize}
\item $\bot^{\pmb{\tau}}:=\mathbf{0}$,
\item $(\phi\vee\psi)^{\pmb{\tau}}:=\phi^{\pmb{\tau}}\otimes\psi^{\pmb{\tau}}$,
\item for the other cases, $\phi^{\pmb{\tau}}:=\phi^\tau$ as in Definition \ref{trans_tipl2ik}.
\end{itemize}
\end{definition}

Note that under such defined translation, the negation $\neg\phi(=\phi\to\bot)$ of \TIPC is translated into \IKt as $\phi^{\pmb{\tau}}\to\mathbf{0}$, which is in general different from the real negation $\phi^{\pmb{\tau}}\to\bot$. We now show that the translation preserves truth over full powerset models and arbitrary teams.

\begin{proposition}\label{tensor_comodel_model}
Let $\MM=(W,R,V)$ be an intuitionistic Kripke model. For any formula $\phi$ and any team $t\subseteq W$, $\MM,t\models \phi\iff \MM^\bullet,t\Vdash\phi^{\pmb{\tau}}$.
\end{proposition}

\begin{proof}
We prove the lemma by induction on $\phi$. 
If $\phi=\bot$, then for any  team $t$, we have $\MM,t\models\bot$ iff $t=\varnothing$ iff $\MM^\bullet,t\Vdash\mathbf{0}$.

The case $\phi=\psi\vee\chi$ follows easily from induction hypothesis. The other cases are proved analogously to the corresponding cases in Proposition \ref{comodel_model}.
\end{proof}


\subsection{Generalized team intuitionistic Kripke semantics
}\label{sec:gtipl}

In this section, we introduce general intuitionistic team Kripke models, which are abstractions of the full powerset models defined in the previous section. We then define generalized team semantics for \TIPC-formulas over these models.
The \textsf{Split} axiom will not any more be valid under the new semantics. Two other axioms  of \TIPC, monotonicity and (weak) elimination axiom of the local disjunction $\vee$ (axioms (2) and (3)), turn out to be not valid either. We give counter-examples to illustrate the failure of these three axioms in this section, and we discuss frame conditions validating these axioms in the next section.



Let us now define general intuitionistic team Kripke frames $\mathfrak{F}=(W,R,\succcurlyeq,\Cup)$, which have a partial order $\succcurlyeq$ and a join operation $\Cup$ on $\wp(W)$  resembling and generalizing the superset relation $\supseteq$ and the set-theoretic union $\cup$.

\begin{definition}\label{g_frame_df}
A {\em general intuitionstic team Kripke frame} is a tuple $\mathfrak{F}=(W,R,\succcurlyeq,\Cup)$, where 
\begin{itemize}
\item $(W,R)$ is an intuitionistic Kripke frame,
\item $\Cup$ is a binary operation on $\wp(W)$ such that $(\wp(W),\Cup,\varnothing)$ forms a bounded join-semilattice, and  $\preccurlyeq$ is the induced partial order;
\item and $\succcurlyeq,\Cup$ and the lifting $R^\circ$ of $R$ on $\wp(W)$ 
satisfy the following conditions:
\begin{enumerate}[(a)]
\item $R^\circ$ and $\succcurlyeq$ satisfy (F2),


\item\label{cond_f3like} $t\succcurlyeq\circ R^\circ s$ implies $R[t]\succcurlyeq s$ for any $t,s\subseteq W$,
\item\label{cond_cup} $R[t]\Cup R[s]=R[t\Cup s]$  for any $t,s\subseteq W$.
\end{enumerate}
\end{itemize}

A {\em general intuitionistic team Kripke model} $\mathfrak{M}=(W,R,\succcurlyeq,\Cup,V)$ is a general intuitionistic team Kripke frame with a persistent  valuation $V:\Prop\to \wp(W)$ on  $(W,R)$, i.e., $(W,R,V)$ is an intuitionistic Kripke model. 
\end{definition}

It follows from the same argument as in Fact \ref{lifting_R_preorder_M0}, the lifting $R^\circ$ is a pre-order, and if $\MM$ is finite,  $R^\circ$ is a partial order. As before, we sometimes abuse the notation and write simply $R$  for the lifted relation $R^\circ$. It is easy to verify that composition relation $\succcurlyeq\circ R^\circ$ is a pre-order as well (e.g.,  transitivity follows from (F2)).



Conditions (a)-(c) in the above definition are imposed to guarantee the persistence of  \TIPC-formulas under the new semantics (as we will verify in Lemma \ref{more_general_persistence}). In Condition (a) we do not require $R^\circ$ and $\succcurlyeq$ to satisfy the (F1) condition, which was however also imposed for bi-relation intuitionistic Kripke models (see Definition \ref{bi-relation_ml}).  The condition (F1) guarantees the persistence of \IK-formulas $\Diamond\phi$ with the diamond modality $\Diamond$. Since the diamond modality does not play a role in our discussion (see the two translations $(\cdot)^{\tau}$ and $(\cdot)^{\pmb{\tau}}$ in Definitions \ref{trans_tipl2ik} and \ref{trans_bold_tipl2ik}), condition (F1) is then not imposed in our models. Every general intuitionistic team Kripke model $\MM=(W,R,\succcurlyeq,\Cup,V)$ is naturally associated with a model $\MM^\bullet=(\wp(W),R^\circ,\succcurlyeq,\Cup,\varnothing,V^\circ)$ whose  valuation $V^\circ$ is defined as 
\(t\in V^\circ(p)\text{ iff }t\preccurlyeq V(p).\) Clearly, this model $\MM^\bullet$  satisfies all the requirements for a tensored bi-relation intuitionistic Kripke model except for (F1)

Condition  (\ref{cond_f3like}) is a special case of the following weaker condition, introduced by Bo\v{z}i\'{c} and Do\v{s}en \cite{BozicDosen84}, on the  model $\MM^\bullet$: 
\begin{itemize}
\item If $t\succcurlyeq r$ and $rR^\circ s$, then there exists $u\in \wp(W)$ such that $tR^\circ u$ and $u\succcurlyeq s$ (i.e., $\succcurlyeq\circ R^\circ\subseteq R^\circ\circ \succcurlyeq$); see also Figure \ref{Fig_F3}.
\end{itemize}
Condition (b) has an interesting special case: If $t\succcurlyeq s$, then we have $R[t]\succcurlyeq R[s]$, since $t\succcurlyeq sR^\circ R[s]$.
This observation implies that the direction $R[t\Cup s]\succcurlyeq R[t]\Cup R[s]$ of Condition (\ref{cond_cup}) is actually always satisfied, as $t\Cup s\succcurlyeq t,s$.


 \begin{figure}[t]
 \begin{center}
\begin{tikzpicture}[scale=.8, transform shape] 
\node  (1)  at (0, 0) {\normalsize$t$};
\node  (3) at +(0,2) {\normalsize$r$};
\node (4) at +(2,2) {\normalsize$s$};
\node  (2) at +(2,0) {\normalsize$u$};
\draw[line width=1pt, dashed] (2) -- (4);
\draw[line width=1pt] (1) -- (3);
\draw [ ->] (3) -- (4);
\draw [dashed,->] (1) -- (2);


\end{tikzpicture}
 \caption{
 }\label{Fig_F3}
\end{center}
  \end{figure}

If the partial order $\succcurlyeq$ is taken to be the superset relation $\supseteq$ (or the binary operator $\Cup$ is taken to be the set-theoretic union $\cup$), Conditions (a)-(c) are automatically satisfied regardless of how $R$ is defined.
Similarly, if $R$ is taken to be the identity relation $\mathsf{id}=\{(w,w)\mid w\in W\}$ (i.e., the underlying Kripke frame $(W,R)$ is a classical frame), then Conditions (a)-(c) are trivially satisfied regardless of how $\succcurlyeq$ and $\Cup$ are defined. We call such frames {\em classical frames}. Note that classical-based logics with team semantics, such as propositional dependence logic \cite{VY_PD} and inquisitive logic \cite{InquiLog}, are defined essentially over  frames of the form $\mathfrak{F}=(W,\mathsf{id},\supseteq,\cup)$, which are very specific types of general intuitionistic team Kripke frames with classical underlying Kripke frames. 
Below we give a nontrivial example of a general intuitionistic team Kripke frame in which $\succcurlyeq\neq \supseteq$ and $R\neq \mathsf{id}$.

\begin{example}\label{example_gframe}
{\em Consider the frame $\mathfrak{F}=(W,R,\succcurlyeq,\Cup)$, where the underlying intuitionistic Kripke frame $(W,R)$ with $W=\{w,u,v\}$ is depicted as the left figure below, and the bounded semilattice $(\wp(W),\Cup,\varnothing)$  is depicted (upside down) as the right figure below (where the dashed arrows represent the lifting $R^\circ$ with reflexive and transitive arrows omitted).
 
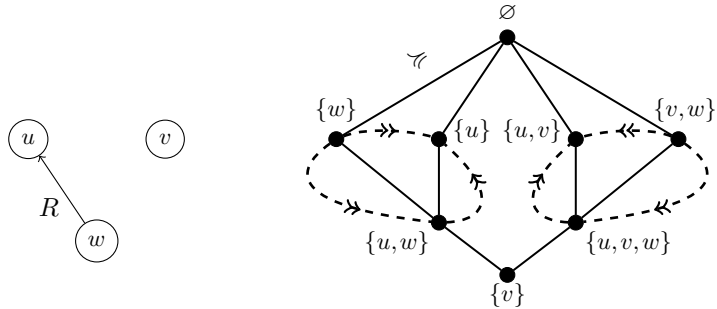
\begin{center}
\begin{tikzpicture}[scale=.9, transform shape]

\node[draw, circle, scale=0.9]  (w0)  at (-6, 0) {\normalsize$w$};
\node[draw, circle, scale=0.9]  (u0) at (-7,1.5) {\normalsize$u$};
\node[draw, circle, scale=0.9] (v0) at (-5,1.5) {\normalsize$v$};
\draw [ ->] (w0) -- (u0);
\node at (-6.7, 0.5) {$R$};

\node at (-1.3,2.7) {\rotatebox{30}{$\succcurlyeq$}};

  \node[draw, circle, scale=0.6, fill=black, label={[below, yshift=-0.1cm]\small $\{v\}$}] (v) at (0,-0.5) {};

  \node[draw, circle, scale=0.6, fill=black, label={[left, xshift=-0.07cm]\small $\{u,v\}$}] (uv) at (1,1.5) {};
 
  \node[draw, circle, scale=0.6, fill=black, label={[right, xshift=0.07cm]\small $\{u\}$}] (u) at (-1,1.5) {};
   \node[draw, circle, scale=0.6, fill=black, label={\small $\{w\}$}] (w) at (-2.5,1.5) {};

   \node[draw, circle, scale=0.6, fill=black, label={[below left, yshift=-0.1cm]\small $\{u,w\}$}] (uw) at (-1,0.27) {};    
   
      \node[draw, circle, scale=0.6, fill=black, label={[xshift=0.1cm]\small $\{v,w\}$}] (vw) at (2.5,1.5) {};
         \node[draw, circle, scale=0.6, fill=black, label={[below right, yshift=-0.1cm]\small $\{u,v,w\}$}] (uvw) at (1,0.27) {};
         
  \node[draw, circle, scale=0.6, fill=black, label={\small $\varnothing$}] (0) at (0,3) {};
         
\draw[line width=.75pt] (v) -- (uw);
\draw[line width=.75pt] (uv) -- (uvw);
\draw[line width=.75pt] (v) -- (uvw);
\draw[line width=.75pt] (uw) -- (w);
\draw[line width=.75pt] (uw) -- (u);
\draw[line width=.75pt] (w) -- (0);
\draw[line width=.75pt] (uv) -- (0);
\draw[line width=.75pt] (u) -- (0);
\draw[line width=.75pt] (vw) -- (0);
\draw[line width=.75pt] (uvw) -- (vw);

\begin{scope}[very thick,decoration={
    markings,
    mark=at position 0.6 with {\arrow{>>}}}
    ] 
\draw[dashed, line width=1pt, postaction={decorate}] (w) .. controls (-3.2,1) and (-3.1,0.5) .. (uw);
\draw[dashed, line width=1pt, , postaction={decorate}]  (w) .. controls (-2,1.7) and (-1.3,1.6) .. (u);
\draw[dashed, line width=1pt, , postaction={decorate}] (vw) .. controls (2,1.7) and (1.3,1.6) .. (uv);
\draw[dashed, line width=1pt, , postaction={decorate}] (vw) .. controls (3.2,1) and (3.1,0.5) .. (uvw);
\end{scope}

\begin{scope}[very thick,decoration={
    markings,
    mark=at position 0.65 with {\arrow{>>}}}
    ] 
\draw[dashed, line width=1pt, , postaction={decorate}] (uw) .. controls (-.3,.3) and (-.1,.7) .. (u);
\draw[dashed, line width=1pt, , postaction={decorate}] (uvw) .. controls (.3,.3) and (.1,0.7) .. (uv);

\end{scope}

  \end{tikzpicture}
\end{center}
It can be verified that $\mathfrak{F}$ indeed forms a general intuitionistic team Kripke frame. We only illustrate the verification of Conditions (a)-(c) for some nontrivial cases; the remaining cases are left to the reader. Note that the full proof contains many symmetric cases, as the above figure of the semilattice $(\wp(W),\Cup,\varnothing)$ is symmetrical vertically, and the teams in the symmetric positions have similar $R^\circ$-successors. 

For Condition (a), assuming $tR r\succcurlyeq s$ we find $t'$ such that $t\succcurlyeq t' R s$. We only consider the case $t=\{w\}$ and $r=\{u,w\}$. Then we have that either  $s=\varnothing$ or $s\in\{\{u,w\},\{w\},\{u\}\}$. In the first case,  take $t'=\varnothing=s$ and we have $t\succcurlyeq \varnothing R\varnothing$. In the second case, take $t'=t=\{w\}$ and we have $t \succcurlyeq tRs$.

For Condition (b), assuming $t\succcurlyeq \circ Rs$ we show that $R[t]\succcurlyeq s$. We only consider the case $t=\{w\}$. For any $s$ such that $t\succcurlyeq \circ R s$, it is easy to observe that  $\{u,w\}\succcurlyeq s$. Thus, $R[\{w\}]=\{u,w\}\succcurlyeq s$.

For Condition (c), we only verify $R[t]\Cup R[s]=R[t\Cup s]$ for any $t=\{w\}$ and $s\neq\varnothing,t$. We distinguish three cases according to the value of $s$. Case 1: $s=\{v\}$: Then 
\[R[\{w\}]\Cup R[\{v\}]=\{u,w\}\Cup \{v\}=\{v\}=R[\{v\}]=R[\{w\}\Cup \{v\}].\]
Case 2: $s\in \{\{u\},\{u,w\}\}$. Since $R[s]=s$ and $R[\{u,w\}]=\{u,w\}$, we have
\[R[\{w\}]\Cup R[s]=\{u,w\}\Cup s=\{u,w\}=R[\{u,w\}]=R[\{w\}\Cup s].\]
Case 3: $s\in  \{\{u,v\},\{v,w\},\{u,v,w\}\}$: Since $\varnothing\neq R[s]\succcurlyeq \{u,v,w\}$ and $\{w\}\Cup s=\{v\}$,  we have
\[R[\{w\}]\Cup R[s]=\{u,w\}\Cup R[s]=\{v\}=R[\{v\}]=R[\{w\}\Cup s].\]
This completes our illustration of the proof. }
\end{example}

We now define the generalized team semantics of \TIPC-formulas over general intuitionistic team Kripke models. The definition is almost the  same as in Definition \ref{df_ts_standard}, except that we now have the more general partial order $\succcurlyeq$ and join $\Cup$ in place of the superset relation $\supseteq$ and the set-theoretic union $\cup$, respectively. 

\begin{definition}\label{gen_team_sem_df}
Let $\mathfrak{M}=(W , R, \succcurlyeq, \Cup,V)$ be a general intuitionistic team Kripke model and $t\subseteq W$ a team. The satisfaction relation $\mathfrak{M},t\models\phi$ is defined inductively exactly as in Definition \ref{df_ts_standard} except for the following cases:
\begin{itemize}
\item $\mathfrak{M},t\models p$ ~~iff~~ $t\preccurlyeq V(p)$ 
\item $\mathfrak{M},t\models \phi\vee\psi$ ~~iff~~ there are $s,r\subseteq W$ such that $t\preccurlyeq s\Cup r$, $\mathfrak{M},s\models \phi$ and $\mathfrak{M},r\models \psi$
\item $\mathfrak{M},t\models \phi\to\psi$ ~~iff~~ for all $s\subseteq W$ with $t \succcurlyeq\circ R s$, $\mathfrak{M},s\models \phi\Longrightarrow\mathfrak{M},s\models \psi$
\end{itemize}
We write $\Gamma\models^{\mathsf{g}}\phi$ (or simply $\Gamma\models\phi$) if for all  general intuitionistic team Kripke models $\mathfrak{M}$ and all teams $t$, $\mathfrak{M},t\models \psi$ for all $\psi\in\Gamma$ implies $\mathfrak{M},t\models \phi$. 
\end{definition}

Let us remark that in the setting of our generalized team semantics, the local disjunction $\vee$ can be understood as a binary diamond modality, for which the corresponding ternary accessibility relation $R_\vee$ is defined as 
\[R_\vee(t,s,r)\iff t\preccurlyeq s\Cup r.\] 

We now verify that \TIPC-formulas are persistent under the new semantics.

\begin{lemma}[Persistence]\label{more_general_persistence}
Let $\mathfrak{M}$ be a model, $t$ and $s$ teams, and $\phi$ a formula.
 If $\mathfrak{M},t\models\phi$ and $t\succcurlyeq \circ Rs$, then $\mathfrak{M},s\models \phi$. In particular,  if $\mathfrak{M},t\models\phi$ and $t\succcurlyeq s$, then $\mathfrak{M},s\models \phi$.
\end{lemma}
\begin{proof}
We prove the lemma by induction on $\phi$. If $\phi=\bot$ and $t\models\bot$, then  $t=\varnothing$. For any $s$ such that $\varnothing\succcurlyeq\circ Rs$, since $\varnothing$ is the greatest element for the order $\succcurlyeq$, we must have that $s=\varnothing$. Thus, $s\models\bot$, as required.

If $\phi=p$, suppose $t\models p$ and $t\succcurlyeq r Rs$. Then $V(p)\succcurlyeq t\succcurlyeq rRs$, which by Condition (\ref{cond_f3like}) and persistence of $V$ implies that $V(p)=R[V(p)]\succcurlyeq s$. Thus, $s\models p$.

If $\phi=\psi\to\chi$, suppose $\mathfrak{M},t\models\psi\to\chi$ and $t\succcurlyeq rRs$. Suppose also $s\succcurlyeq uRv$ and $v\models\psi$.  By (F2), there exists $r'\subseteq W$ such that $r\succcurlyeq r'Ru$. Thus, $t\succcurlyeq  r\succcurlyeq r'RuRv$, which by assumption implies that $v\models\chi$.

%


If $\phi=\psi\vee\chi$, suppose $\mathfrak{M},t\models\psi\vee\chi$ and $t\succcurlyeq r Rs$. Then there exist $t_1,t_2\subseteq W$ such that $t\preccurlyeq t_1\Cup t_2$, $t_1\models\psi$ and $t_2\models\chi$.  Since $t_1\Cup t_2\succcurlyeq t\succcurlyeq rRs$, by Conditions (\ref{cond_f3like}) and (\ref{cond_cup}), we have that $R[t_1]\Cup R[t_2]=R[t_1\Cup t_2]\succcurlyeq s$. Moreover, by definition we have 
$t_1\succcurlyeq t_1RR[t_1]$ and $t_2\succcurlyeq t_2RR[t_2]$. 
Thus, by induction hypothesis, we have that $R[t_1]\models\psi$ and $R[t_2]\models\chi$, from which it follows that  $s\models \psi\vee\chi$ follows.
%
%
%

The other cases follow immediately from the induction hypothesis.
\end{proof}

The deduction theorem is an immediate corollary of persistence.

\begin{corollary}[Deduction theorem]
For any set $\Gamma\cup\{\phi,\psi\}$ of formulas, we have that
\[\Gamma,\phi\models^{\textsf{g}}\psi\iff \Gamma\models^{\textsf{g}}\phi\to\psi.\]
\end{corollary}

%
%
%

Next, we prove that the empty team property holds for \TIPC-formulas under the generalized semantics.

\begin{lemma}[Empty team property]\label{gTIPL_emptyteam_prop}
For any model $\mathfrak{M}$ and formula $\phi$, we have that $\mathfrak{M},\varnothing\models\phi$. 
\end{lemma}
\begin{proof}
The lemma is proved by induction. If $\phi=\bot$, then $\varnothing\models\bot$ by definition. If $\phi=p$, then $\varnothing\models p$ follows from the fact that $\varnothing\preccurlyeq V(p)$ (since $\varnothing$ is the greatest element  with respect to $\succcurlyeq$). 

If $\phi=\psi\vee\chi$, by induction hypothesis we have $\varnothing\models\psi$ and $\varnothing\models\chi$. Since $\varnothing$ is the identity of the join $\Cup$, we have $\varnothing\preccurlyeq\varnothing=\varnothing\Cup\varnothing$. Hence $\varnothing\models\psi\vee\chi$. 

If $\phi=\psi\to\chi$, for any $s,r$ such that $\varnothing\succcurlyeq r R s$ and $s\models\psi$, 
we must have $r=\varnothing$ and thus $s=\varnothing$. Hence, by induction hypothesis, we conclude $s\models\chi$.

The other cases follow immediately from induction hypothesis.
\end{proof}


The persistence property justifies that our generalized semantics is a well-defined intuitionistic semantics. Together with the empty team property these suggest that our generalized semantics inherits certain nice features of the standard team semantics of \TIPC. 
Let us now examine the logic that the generalized semantics induces.

It is easy to verify that  all \IPC axioms are sound for the language $[\bot,\wedge,\vvee,\to]$ 
over  general intuitionistic team Kripke models (since, in particular,  persistence, (F2), and the empty team property hold). Axioms (1) and (4)-(6) of \TIPC (see Definition \ref{tipc_df}) also remain sound with respect to the generalized team semantics, as we verify in the following.

\begin{lemma}\label{gener_fact_vvee}
Let $\phi,\psi$ and $\chi$ be formulas.
\begin{enumerate}[(i)]
\item\label{gener_fact_vvee_semilattice}  $\phi\vee\psi\equiv^{\textsf{g}}\psi\vee\phi$ and $(\phi\vee\psi)\vee\chi\equiv^{\textsf{g}}\phi\vee(\psi\vee\chi)$.
\item\label{gener_fact_vvee2vee} $\phi\models^{\textsf{g}}\phi\vee\psi$ and $\phi\vvee\psi\models^{\textsf{g}}\phi\vee\psi$.
\item\label{gener_fact_dstr} $\phi\vee(\psi\vvee\chi)\equiv^{\textsf{g}}(\phi\vee\psi)\vvee(\phi\vee\chi)$.
\end{enumerate}
\end{lemma}
\begin{proof}
The first clause of item (\ref{gener_fact_vvee_semilattice}) follows directly from the commutativity of $\Cup$. For the second clause, we only show $(\phi\vee\psi)\vee\chi\models^{\textsf{g}}\phi\vee(\psi\vee\chi)$; the other direction is similar. Suppose $\MM,t\models (\phi\vee\psi)\vee\chi$ for some $\MM=(W,R,\succcurlyeq, \Cup,V)$. Then there exist $r,s\subseteq W$ such that $t\preccurlyeq r\Cup s$, $r\models\phi\vee\psi$ and $s\models\chi$. Since $r\models\phi\vee\psi$, there exist $r_0,r_1\subseteq W$ such that $r\preccurlyeq r_0\Cup r_1$, $r_0\models\phi$ and $r_1\models \psi$. Thus, we have $r_0\Cup (r_1\Cup s)\models\phi\vee(\psi\vee\chi)$. To conclude that $t\models \phi\vee(\psi\vee\chi)$, by persistence, it suffices to show that $t\preccurlyeq r_0\Cup (r_1\Cup s)$, which reduces to showing $t\preccurlyeq (r_0\Cup r_1)\Cup s$, as $\Cup$ is associative. Now, since $r\preccurlyeq r_0\Cup r_1$ and $\preccurlyeq $ is the order of the semilattice $(\wp(W),\Cup,\varnothing)$, we obtain $t\preccurlyeq r\Cup s\preccurlyeq (r_0\Cup r_1)\Cup s$, as required.

The second clause of item (\ref{gener_fact_vvee2vee}) is a consequence of the first clause. For the first clause of item (\ref{gener_fact_vvee2vee}), suppose $t\models \phi$. By the empty team property (Lemma \ref{gTIPL_emptyteam_prop}), we have $\varnothing\models\psi$. Since $\MM,t\preccurlyeq t=t\Cup \varnothing$, we conclude that $t\models\phi\vee\psi$. 

Item (\ref{gener_fact_dstr}) follows directly from definition.
\end{proof}

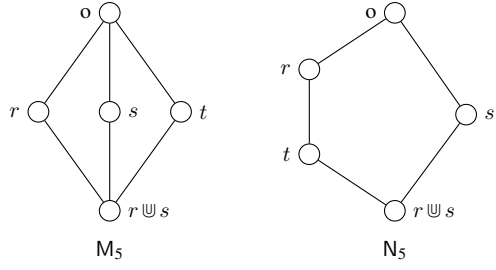
\begin{figure}[t]
 \begin{center}
\begin{tikzpicture}[scale=.75, transform shape] 
\node[draw, circle, label=left: $\omicron$]  (1)  at (0, 0) {};
\node[draw,  circle, label=left: $r$]  (2) at (-1.25,-1.75) {};
\node[draw,  circle, label=right: $s$] (3) at (0,-1.75) {};
\node[draw, circle, label=right: $t$]   (4) at (1.25,-1.75) {};
\node[draw, circle, label=right: $r\Cup s$]  (6) at (0,-3.5) {};
\draw (1) -- (2);
\draw (1) -- (3);
\draw (1) -- (4);
\draw (2) -- (6);
\draw (3) -- (6);
\draw (4) -- (6);

\node at (0,-4.2) {$\mathsf{M}_5$};

\node[draw, circle, label=left: $\omicron$]  (1a)  at (5, 0) {};
\node[draw, circle, label=right: $s$]  (2a) at (6.25,-1.8) {};
\node[draw, circle, label=left: $r$] (3a) at (3.5,-1) {};
\node[draw, circle, label=right: $r\Cup s$]   (4a) at (5,-3.5) {};
\node[draw, circle, label=left: $t$]  (6a) at (3.5,-2.5) {};
\draw (1a) -- (2a);
\draw (1a) -- (3a);
\draw (2a) -- (4a);
\draw (3a) -- (6a);
\draw (4a) -- (6a);

\node at (5,-4.2) {$\mathsf{N}_5$};

\end{tikzpicture}
 \caption{Lattices $\mathsf{M}_5$ and $\mathsf{N}_5$ (placed upside down)
 }\label{M5N5ctexam}
\end{center}
\end{figure}

However, as expected, not all axioms of \TIPC  are sound in the generalized team semantics. We now give counter-examples to illustrate  the failure of the three remaining axioms of \TIPC: Axioms (2) and (3), and the \textsf{Split} axiom. For each axiom we provide two  counter-models, whose underlying frames are denoted as $\mathfrak{F}_{\mathsf{M}}$ and $\mathfrak{F}_{\mathsf{N}}$. The frame $\mathfrak{F}_{\mathsf{M}}=(W,R,\succcurlyeq,\Cup)$ has as domain  an arbitrary finite set $W$ with at least 3 elements, and the partial order $R=\textsf{id}$ (i.e., the underlying Kripke frame $(W,R)$ is classical, and thus Conditions (a)-(c) in \Cref{g_frame_df} are trivially satisfied). The associated (finite) bounded semilattice $(\wp(W),\Cup,\varnothing)$ of $\mathfrak{F}_{\mathsf{M}}$ is thus  a lattice (with the meet $\Cap$ defined as $t\Cap s=\bigCup\{r\in W\mid r\preccurlyeq t,s\}$), and we require that it contains the lattice $\mathsf{M}_5$ as shown (upside down) in \Cref{M5N5ctexam} as a sublattice. The finite frame $\mathfrak{F}_{\mathsf{N}}=(W,R,\succcurlyeq,\Cup)$ is defined the same way as $\mathfrak{F}_{\mathsf{M}}$, except that its  bounded (semi)lattice $(\wp(W),\Cup,\varnothing)$ contains instead the lattice $\mathsf{N}_5$ as shown (upside down) in \Cref{M5N5ctexam} as a sublattice.  Recall (see, e.g., \cite{BurrisSankappanavar2012}) that a lattice is distributive if and only if none of its sublattices is isomorphic to $\mathsf{M}_5$ or $\mathsf{N}_5$ (note that $\mathsf{M}_5$ and $\mathsf{N}_5$ are horizontally symmetric, and thus the direction of the order $\succcurlyeq$ in \Cref{M5N5ctexam} does not make a substantial difference in the discussion). Thus, the lattices $(\wp(W),\Cup,\Cap,\varnothing)$ for $\mathfrak{F}_{\mathsf{M}}$ and $\mathfrak{F}_{\mathsf{N}}$ are clearly not distributive. The reason for the slightly more involved construction of these counter-models will become clear in the next section, when we discuss distributive frames. Similar counter-examples to our \Cref{ctexam_failureCup_closure} below can also be found in \cite{Puncochar17} in a different and algebraic setting.

We now first give our two counter-examples for Axiom (3)  $(\phi\to\chi)\to(\phi\vee\psi\to\chi\vee\psi)$, which expresses the monotonicity of the local disjunction $\vee$. 

\begin{ctexample}[Axiom (3)]\label{ctexam_axm3}
 {\em (1). Consider the general intuitionistic team Kripke model $\MM=(\mathfrak{F}_{\mathsf{M}},V)$, where the frame $\mathfrak{F}_{\mathsf{M}}$ is as described above, the valuation $V$ satisfies
\[V(p)=r, ~V(q)=s\text{ and }V(p')=\omicron,\] 
and the teams $\omicron, r,s,t,r\Cup s$  are as shown in $\mathsf{M}_5$ in \Cref{M5N5ctexam}.
Note that since $R=\mathsf{id}$, the valuation $V$ is trivially persistent.
Moreover, since $\omicron\models p,q,p'$, the valuation is also legitimate (or consistent with the empty team property) in case $\omicron=\varnothing$.

We have that $\MM,t\not\models (p\to p')\to (p\vee q\to p'\vee q)$. To see why, first observe that since $R=\textsf{id}$,  a $\succcurlyeq\circ R$-successor of any point is a $\succcurlyeq$-successor, and vice versa. We first show that $t\models p\to p'$. Suppose $t\succcurlyeq t'$ and $t'\models p$. By definition, we have $t'\preccurlyeq V(p)=r$. Since $\mathsf{M}_5$ (with the meet $\Cap_{\mathsf{M}_5}$) is a sublattice of $(\wp(W),\Cup,\Cap,\varnothing)$,  $\omicron=r\Cap_{\mathsf{M}_5} t=r\Cap t\succcurlyeq t'$. Now, since $V(p')=\omicron$, we obtain $t'\models p'$, as required.

Next,  we have $t\models p\vee q$, since  $r\models p$, $s\models q$, and $t\preccurlyeq r\Cup s$. Finally, we verify that $t\not\models p'\vee q$, from which we would conclude that $t\not\models p\vee q\to p'\vee q$. For any $r',s'$ such that $r'\models p'$, $s'\models q$, we have, by definition, that $r'\preccurlyeq V(p')=\omicron\preccurlyeq s$ and $s'\preccurlyeq V(q)=s$, which imply $s\succcurlyeq r'\Cup s'$. Thus, we have that $t\not\preccurlyeq r'\Cup s'$, since otherwise we would have $t \preccurlyeq s$, which is not the case.
}

\vspace{\baselineskip}

\noindent  {\em (2). Consider the general intuitionistic team Kripke model $\MM=(\mathfrak{F}_{\mathsf{N}},V)$, where  the frame $\mathfrak{F}_{\mathsf{N}}$ is as described above,   the valuation $V$ satisfies
\[V(p)=s,~V(p')=\omicron\text{ and } V(q)=r,\] 
and the teams $\omicron, r,s,t,r\Cup s$  are as shown in $\mathsf{N}_5$ in \Cref{M5N5ctexam}.

We have that $\MM,t\not\models (p\to p')\to (p\vee q\to p'\vee q)$. To see why, we first show that $t\models p\to p'$. For any $t'\preccurlyeq t$ with $t'\models p$, we have $t'\preccurlyeq V(p)=s$. But since $\mathsf{N}_5$ is a sublattice of $(\wp(W),\Cup,\Cap,\varnothing)$, we have that $t'\preccurlyeq t\Cap s=t\Cap_{\mathsf{N}_5}s=\omicron=V(p')$. Thus, $t'\models p'$.

Next, clearly $t\models p\vee q$, since $t\preccurlyeq r\Cup s$, $r\models q$  and $s\models p$. Finally, we verify that $t\not\models p'\vee q$, from which it would follow that $t\not\models p\vee q\to p'\vee q$. For any $r',s'$ such that  $r'\models q$ and $s'\models p'$, we have, by definition, that $r'\preccurlyeq V(q)=r$ and $s'\preccurlyeq V(p')=\omicron\preccurlyeq r$, which imply that $r'\Cup s'\preccurlyeq r\preccurlyeq t$. Since $r\neq t$, we obtain $t\not\preccurlyeq r'\Cup s'$.
}
\end{ctexample}




The above counter-examples show the failure of the monotonicity of the local disjunction $\vee$, stated as an implication $(\phi\to\chi)\to(\phi\vee\psi\to\chi\vee\psi)$. Nevertheless, it is straightforward to verify that the following (weaker form of) monotonicity of $\vee$ does hold: 
\[\text{If }\phi\models^{\textsf{g}}\chi,\text{ then }\phi\vee\psi\models^{\textsf{g}}\chi\vee\psi.\]
In other words, the following deduction rule of $\vee$ is sound:
\[
\AxiomC{$\phi\vdash\chi$}\RightLabel{$\vee\textsf{Mon}$} \UnaryInfC{$\phi\vee\psi\vdash\chi\vee\psi$}\DisplayProof
\]
Note that we do not allow any context set $\Gamma$ to occur on the left-hand side of $\vdash$ in the above rule.


Next, we shall illustrate the failure of Axiom (2) $(\phi\to\alpha)\to ((\psi\to\alpha)\to(\phi\vee\psi\to\alpha))$. We do so by first giving two examples to illustrate the failure of closure under $\Cup$-join for standard formulas (i.e., $\vvee$-free formulas). Recall that in the standard team semantics, standard formulas are, however, closed under $\cup$-join (or unions). We say that a formula $\phi$ is {\em closed under (finite) join} if for any model $\MM=(W,R,\succcurlyeq,\Cup,V)$ and teams $t,s\subseteq W$
\[\MM,t\models\phi\text{ and }\MM,s\models\phi\Longrightarrow \MM,t\Cup s\models\phi.\tag{\textsf{Join Closure Property}}\]


\begin{ctexample}[Failure of join closure for standard formulas]\label{ctexam_failureCup_closure} 
{\em (1). Consider the general intuitionistic team Kripke model $\MM=(\mathfrak{F}_{\mathsf{M}},V)$, where the valuation $V$ satisfies
\[V(p)=t\text{ and }V(q)=\omicron.\]
Clearly, $r\Cup s\not\models p\to q$, since for the team $t$, we have $r\Cup s\succcurlyeq t\models p$, whereas $t\not\models q$. On the other hand, we have $r\models p\to q$ and $s\models p\to q$. 

We now give the proof of $r\models p\to q$; the proof of $s\models p\to q$ is similar. Suppose $r'\preccurlyeq r$ and $r'\models p$. Then, by definition, we have that $r'\preccurlyeq V(p)=t$. Since $\mathsf{M}_5$ is a sublattices of $(\wp(W),\Cup,\Cap,\varnothing)$, we have that $r'\preccurlyeq r\Cap t=r\Cap_{\mathsf{M}_5}t=\omicron=V(q)$, which implies $r'\models q$, as required.

\vspace{\baselineskip}

\noindent (2). Consider the general intuitionistic team Kripke model $\MM=(\mathfrak{F}_{\mathsf{N}},V)$, where  the valuation $V$ satisfies
\[V(p)=t\text{ and }V(q)=r.\]
Clearly, $r\Cup s\not\models p\to q$, for a similar reason as above. On the other hand, we have $r\models p\to q$ and $s\models p\to q$. The former follows easily from persistence and the fact that $r\models p$ and $r\models q$. The latter is proved by going through a similar argument to that in (1) where the fact $t\Cap s=\omicron$ is used.
}

\end{ctexample}

Closure under (finite) joins of a formula $\phi$ in a model is actually equivalent to validity of the {\em idempotent law} $(\phi\vee\phi)\leftrightarrow\phi$ in the model, as we now show. Note that $\phi\to\phi\vee\phi$ is always valid (by \Cref{gener_fact_vvee}(\ref{gener_fact_vvee2vee}) and the deduction theorem).

\begin{lemma}\label{eq_Ucl_idem}
Let $\MM=(W,R,\succcurlyeq,\Cup,V)$ be a general intuitionistic team Kripke model. For any formula $\phi$,
\[\phi\text{ is closed under (finite) joins}\text{ in }\MM\iff \MM\models\phi\vee\phi\to\phi.\]
\end{lemma}
\begin{proof}
``$\Longrightarrow$'': For any team $t\subseteq W$, if $\MM,t\models\phi\vee\phi$, then there exist $r,s\subseteq W$ such that $t\preccurlyeq r\Cup s$, $r\models\phi$ and $s\models\phi$. Since $\phi$ is closed under joins in $\MM$, we have that $r\Cup s\models \phi$, which then implies $t\models\phi$ by persistence. 

``$\Longleftarrow$'': Suppose $\MM,s\models\phi$ and $\MM,t\models\phi$ for some $s,t\subseteq W$. By definition $s\Cup t\models\phi\vee\phi$. Now, since $\MM\models\phi\vee\phi\to\phi$, we conclude that $s\Cup t\models\phi$.
\end{proof}

By \Cref{IPL_flat}, over the standard team semantics the idempotent law $\alpha\vee\alpha\to\alpha$ for standard formulas $\alpha$ is a theorem in \TIPC. By the above lemma, Counter-example \ref{ctexam_failureCup_closure} also shows that $\alpha\vee\alpha\to\alpha$ is no longer sound under the generalized semantics. In particular, in the two models of the counter-example we have $r\Cup s\models (p\to q)\vee(p\to q)$, whereas $r\Cup s\not\models p\to q$. 
The same example also illustrates that the substitution instance of Axiom (2) $(\phi\to\alpha)\to((\phi\to\alpha)\to (\phi\vee\phi\to\alpha))$ obtained by putting $\phi=\alpha=p\to q$ is not sound, since, 
e.g., in both models in Counter-example \ref{ctexam_failureCup_closure}, we have $r\Cup s\not\models (p\to q)\vee(p\to q)\to (p\to q)$. 
From this it also follows that $ (p\to q)\vee(p\to q)\not\models p\to q$, even though $p\to q\models p\to q$. Therefore, unlike with Axiom (3), the (weaker) rule that corresponds to Axiom (2) ``from $\phi\vdash\alpha$ and $\psi\vdash\alpha$, derive $\phi\vee\psi\vdash\alpha$'' is not sound either.

Finally, we also construct two counter-models for the \textsf{Split} axiom, whose underlying frames are two concrete instances of the frames $\mathfrak{F}_{\mathsf{M}}$ and $\mathfrak{F}_{\mathsf{N}}$. We, however, did not succeed to construct these counter-models with arbitrary frames $\mathfrak{F}_{\mathsf{M}}$ and $\mathfrak{F}_{\mathsf{N}}$.

\begin{figure}[t]
 \begin{center}
\begin{tikzpicture}[scale=.9, transform shape] 

\node at (-7.7,2.7) {\rotatebox{30}{$\succcurlyeq$}};

\node[draw, circle, fill=black, label=above: $\varnothing$, scale=0.8]  (1)  at (-7, 3) {};
\node[draw, circle, fill=black, label=left: $r$, scale=0.8]  (2) at (-8.3,1.8) {};
\node[draw, circle, fill=black, label=right: $s$, scale=0.8] (3) at (-7,1.8) {};
\node[draw, circle, fill=black, label=right: $t$, scale=0.8]   (4) at (-5.7,1.8) {};
\node[draw, circle, scale=0.8]  (5) at (-8.3,0.6) {};
\node[draw, circle, fill=black, label=right: $r\Cup s$, scale=0.8]  (6) at (-7,0.6) {};
\node[draw, circle, scale=0.8]  (7) at (-5.7,0.6) {};
\node[draw, circle,  scale=0.8]  (8) at (-7,-.6) {};
\draw (1) -- (2);
\draw (1) -- (3);
\draw (1) -- (4);
\draw (2) -- (5);
\draw (2) -- (6);
\draw (3) -- (6);
\draw (4) -- (6);
\draw (4) -- (7);
\draw (5) -- (8);
\draw (6) -- (8);
\draw (7) -- (8);


\node at (-1.3,2.7) {\rotatebox{30}{$\succcurlyeq$}};

  \node[draw, circle, scale=0.8, fill=black, label=right: $r\Cup s$] (v) at (0,-0.5) {};

  \node[draw, circle, scale=0.8] (uv) at (1,1.5) {};
 
  \node[draw, circle, scale=0.8, fill=black, label={[right, xshift=0.07cm]\small $r$}] (u) at (-1,1.5) {};
   \node[draw, circle, scale=0.8] (w) at (-2.5,1.5) {};

   \node[draw, circle, scale=0.8, fill=black, label={[below left, yshift=-0.1cm]\small $t$}] (uw) at (-1,0.27) {};    
   
      \node[draw, circle, scale=0.8] (vw) at (2.5,1.5) {};
         \node[draw, circle, scale=0.8, fill=black, label={[below right, yshift=-0.1cm]\small $s$}] (uvw) at (1,0.27) {};
         
  \node[draw, circle, scale=0.8, fill=black, label={\small $\varnothing$}] (0) at (0,3) {};
         
\draw (v) -- (uw);
\draw (uv) -- (uvw);
\draw (v) -- (uvw);
\draw (uw) -- (w);
\draw (uw) -- (u);
\draw (w) -- (0);
\draw (uv) -- (0);
\draw (u) -- (0);
\draw (vw) -- (0);
\draw (uvw) -- (vw);

\end{tikzpicture}
 \caption{Two bounded join-semilattices $(\wp(W),\Cup,\varnothing)$ (placed upside down)  with $\mathsf{M}_5$ and $\mathsf{N}_5$ as a sublattice respectively.}\label{Fig_example_nonregular}
\end{center}
\end{figure}
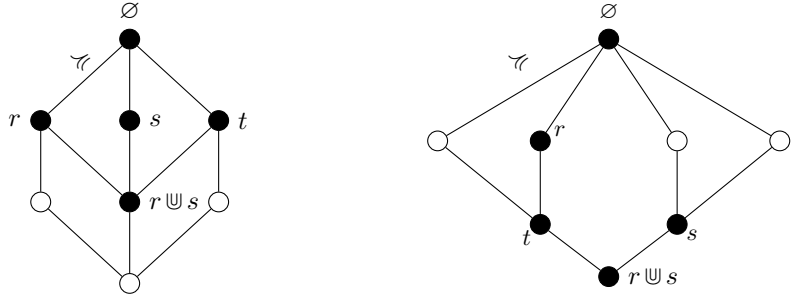

\begin{ctexample}[\textsf{Split}]\label{ctexam_split2} 
{\em (1). Consider the general intuitionistic team Kripke model $\MM=(W,R,\succcurlyeq,\Cup,V)$, where $W$ is an arbitrary $3$-element set,  $R=\textsf{id}$, the bounded semilattice $(\wp(W),\Cup,\varnothing)$ is the left figure shown (upside down) in \Cref{Fig_example_nonregular}, and  the valuation $V$ satisfies
\[V(p)=t,~V(q)=\varnothing,~V(a)=r~\text{ and }~V(b)=s.\]
Note that the highlighted points in $(\wp(W),\Cup,\varnothing)$ form a sublattice, the lattice $\mathsf{M}_5$.
Then $r\Cup s\not\models ((p\to q)\to a\vvee b)\to((p\to q)\to a)\vvee((p\to q)\to b)$.
To see why, observe first that $r\Cup s\models (p\to q)\to a\vvee b$, since the only teams $t' \preccurlyeq r\Cup s$ such that $t'\models p\to q$ are $r,s,$ and $\varnothing$, on which either $a$ or $b$ is true.  But $r\Cup s\not\models (p\to q)\to a$ and $r\Cup s\not\models (p\to q)\to b$, since $s\not\models a$, and  $r\not\models b$.

\vspace{\baselineskip}

\noindent (2). Consider the general intuitionistic team Kripke model $\MM=(W,R,\succcurlyeq,\Cup,V)$, where $W$ is an arbitrary $3$-element set,  $R=\textsf{id}$, the bounded semilattice $(\wp(W),\Cup,\varnothing)$ is the right figure shown (upside down) in \Cref{Fig_example_nonregular}, and  the valuation $V$ satisfies
\[V(p)=t,~V(q)=r=V(a)~\text{ and }~V(b)=s.\]
Note that the semilattice $(\wp(W),\Cup,\varnothing)$ is the same  as the one in Example \ref{example_gframe}, and the highlighted points form a sublattice, the lattice $\mathsf{N}_5$.
By a similar argument to the above, it can be verified that $r\Cup s\not\models ((p\to q)\to a\vvee b)\to((p\to q)\to a)\vvee((p\to q)\to b)$.
}
\end{ctexample}

Consider
 the logic \(\TGIPC:=\{\phi\mid \mathsf{F}^{\mathsf{g}}\models\phi\}\)  of the class $\mathsf{F}^{\mathsf{g}}$ of all general intuitionistic team Kripke frames.
By the observations above, all axioms of \TIPC except for Axioms (2) and (3) and the \textsf{Split} axiom belong to \TGIPC. Whether these axioms actually axiomatize the logic $\TGIPC$ is left as future work.


\subsection{Recovering the standard \TIPC, and distributivity}\label{sec:recover}

In this section, we show that over  general intuitionistic team Kripke frames that are distributive,  Axioms (2) and (3) and the \textsf{Split} axiom are sound, and thus all \TIPC axioms are recovered over such frames; in other words, distributivity is a sufficient condition for validating these axioms. It will then follow also from the counter-examples provided at the end of the previous section that distributivity is actually the necessary and sufficient condition for validating Axioms (2) and (3) over the class of finite classical frames. We also show that over distributive frames standard formulas satisfy the join closure property. 
Similar results are also  found in related work in \cite{Puncochar17} in a different and simpler  setting; we will discuss these connections in more detail in Section \ref{sec:conclusion}.

Let us first define the distributivity condition.

\begin{definition}
A general intuitionistic team Kripke frame $\mathfrak{F}=(W,R,\succcurlyeq,\Cup)$ is called {\em distributive} if the join-semilattice $(\wp(W),\Cup, \varnothing)$ is  distributive, namely 
\[t\preccurlyeq r\Cup s\Longrightarrow \exists r',s'\subseteq W: r'\preccurlyeq r,~ s'\preccurlyeq s,~\text{ and }~t =r'\Cup s'.\]
We denote by $\mathsf{DF}$ the class of distributive general intuitionistic team Kripke frames. 
\end{definition}

A frame $\mathfrak{F}=(W,\supseteq,R,\cup,\varnothing)$ with $R$ a partial order 
is always distributive, as $(\wp(W),\cup,\varnothing)$ is clearly a distributive semilattice. Let us  elaborate on the notion of distributivity. 
Recall, from e.g., \cite{BurrisSankappanavar2012}, that any distributive join-semilattice in which meets exist is a distributive lattice. Meets of finite join-semilattices always exist. A lattice is distributive if and only if none of its sublattices is isomorphic to $\mathsf{M}_5$ or $\mathsf{N}_5$ shown in \Cref{M5N5ctexam} (placed upside down with the highlighted nodes $r,s,t$ witnessing the failure of the distributivity condition). We have essentially already pointed out that the two types of finite frames $\mathfrak{F}_{\mathsf{M}}$ and $\mathfrak{F}_{\mathsf{N}}$ we constructed in the previous section for Counter-examples \ref{ctexam_axm3} and \ref{ctexam_failureCup_closure} (as well as the two concrete frames constructed in Counter-example \ref{ctexam_split2}) are thus not distributive. Furthermore, any finite $k$-element non-distributive classical frame (i.e., $R=\textsf{id}$) is thus isomorphic to either $\mathfrak{F}_{\mathsf{M}}$ or $\mathfrak{F}_{\mathsf{N}}$ with $k$ elements. Our Counter-examples \ref{ctexam_axm3} and \ref{ctexam_failureCup_closure} and the discussions afterwards in the previous section thus show that over finite classical frames, distributivity is a necessary condition for validating Axioms (2) and (3), the idempotent law of the local disjunction $\vee$ for standard formulas, and join closure of standard formulas.

\begin{corollary}\label{distr_nec_cond}
For any finite classical frame $\mathfrak{F}$ that is not in $\mathsf{DF}$, we have that Axiom (2), Axiom (3) and the idempotent law $\alpha\vee\alpha\to\alpha$ (with $\alpha$ standard) are not valid on $\mathfrak{F}$, and standard formulas are not closed under joins over $\mathfrak{F}$.
\end{corollary}

Whether distributivity is a necessary condition for the same axioms for the local disjunction $\vee$ and join closure property also over non-classical frames (i.e., frames in which $R\neq\textsf{id}$) is unclear. As for the \textsf{Split} axiom for the global disjunction $\vvee$, we conjecture that distributivity is not a necessary condition for validating the axiom even over finite classical frames. Our Counter-example \ref{ctexam_split2} from the previous section was built on two concrete non-distributive frames of $\mathfrak{F}_{\mathsf{M}}$ or $\mathfrak{F}_{\mathsf{N}}$ type; it is not clear how to generalize the argument to arbitrary $\mathfrak{F}_{\mathsf{M}}$ or $\mathfrak{F}_{\mathsf{N}}$ type frames.

We now proceed to show that distributivity is, on the other hand, a sufficient condition for validating all above-mentioned axioms and properties, over arbitrary frames. We start with Axioms (3).


\begin{proposition}\label{gener_a3}
For any formulas $\phi,\psi$ and $\chi$, we have that
$\mathsf{DF}\models(\phi\to\chi)\to(\phi\vee\psi\to\chi\vee\psi)$.
\end{proposition}
\begin{proof}
Let $\MM=(W,R,\succcurlyeq,\Cup,V)$ be a model with a distributive underlying frame. Suppose $\MM,t\models \phi\to\chi$ for some $t\subseteq W$. Suppose also that $t\succcurlyeq \circ Rs$ and $s\models \phi\vee\psi$. Then there are $r_0,r_1\subseteq W$ such that $s\preccurlyeq r_0\Cup r_1$, $r_0\models\phi$ and $r_1\models\psi$. By distributivity, there exist $r_0',r_1'\subseteq W$ such that $r_0'\preccurlyeq r_0$, $r_1'\preccurlyeq r_1$, and $s=r_0'\Cup r_1'$. Since $t\succcurlyeq \circ Rs\succcurlyeq  r_0'\preccurlyeq r_0$, by persistence, we have $r_0'\models\phi\to\chi$ and $r_0'\models\phi$, which imply  $r_0'\models\chi$. On the other hand, since $r_1\succcurlyeq r_1'$, we also have $r_1'\models\psi$ by persistence. Hence, we conclude that $s=r_0'\Cup r_1'\models\chi\vee\psi$, as required.
\end{proof}

Next, we show that over  general intuitionistic team Kripke 
models whose underlying frame are distributive, standard formulas  are closed under finite joins in general, and under arbitrary joins given certain additional conditions.

\begin{lemma}[Join closure of standard formulas]\label{Cup_closure_gts}
Let $\MM=(W , R, \succcurlyeq, \Cup,V)$ be a model with a distributive underlying frame.
For any standard formula $\alpha$, any teams $t_0,t_1\subseteq W$,
\[\MM,t_0\models\alpha\text{ and }\MM,t_1\models\alpha\Longrightarrow \MM,t_0\Cup t_1\models\alpha.\]
Moreover, if $(\wp(W),\Cup,\varnothing)$ forms a distributive  complete semilattice, and distributivity and Condition (c) hold for arbitrary $\Cup$-joins, then
\[\MM,t\models\alpha\text{ for every }t\in T\Longrightarrow \MM,\bigsCup T\models\alpha.\]
\end{lemma}
\begin{proof}
We only give the proof for the finite join case. The arbitrary join case follows from essentially the same argument, in which all the steps with finite joins can be easily adapted to the infinite case given the completeness assumption.

We proceed by induction on $\alpha$. If $\alpha=\bot$, then $\mathfrak{M},t_0\models \bot$  and $\mathfrak{M},t_1\models \bot$ imply that $t_0=\varnothing=t_1$. Thus, $t_0\Cup t_1=\varnothing\models\bot$.

If $\alpha=p$, then $t_0,t_1\models p$ implies $t_0,t_1\preccurlyeq V(p)$. Since $(\wp(W),\Cup,\varnothing)$ is a lattice,   $t_0\Cup t_1\preccurlyeq V(p)$ and so  $t_0\Cup t_1\models p$.


Case $\alpha=\beta\vee\gamma$. Suppose $t_0,t_1\models\beta\vee\gamma$. Then there are $r_0,s_0,r_1,s_1\in \wp(W)$ such that $t_0\preccurlyeq r_0\Cup s_0$, $t_1\preccurlyeq r_1\Cup s_1$, $r_0\models\beta$, $s_0\models\gamma$, $r_1\models\beta$, and $s_1\models\gamma$. By induction hypothesis, we have that $r_0\Cup r_1\models\beta$ and $s_0\Cup s_1\models\gamma$. Since $(\wp(W),\Cup,\varnothing)$ is a lattice, 
\[t_0\Cup t_1\preccurlyeq(r_0\Cup s_0)\Cup (r_1\Cup s_1)=(r_0\Cup r_1)\Cup(s_0\Cup s_1)\]
Hence we conclude that $t_0\Cup t_1\models\beta\vee\gamma$.

Case $\alpha=\beta\to\gamma$. Suppose $t_0,t_1\models\beta\to\gamma$. Suppose also that $t_0\Cup t_1\succcurlyeq\circ R s$ and $s\models\beta$. Then, by Conditions (b) and (c), we have $R[t_0]\Cup R[t_1]= R[t_0\Cup t_1]\succcurlyeq s$. By distributivity, there exist $r_0,r_1\subseteq W$ such that $r_0\preccurlyeq R[t_0]$, $r_1\preccurlyeq R[t_1]$,  and $s= r_0\Cup r_1$. Since $t_0R^\circ R[t_0]\succcurlyeq r_0\preccurlyeq s$ and $t_1R^\circ R[t_1]\succcurlyeq r_1\preccurlyeq s$, by  persistence, we have $r_0,r_1\models\beta\to\gamma$ and $r_0,r_1\models\beta$. Thus, $r_0,r_1\models\gamma$.  Hence, by induction hypothesis, we conclude that $s=r_0\Cup r_1\models \gamma$.

The case $\alpha=\beta\wedge\gamma$ is straightforward.
\end{proof}

Let us remark that in the above proof, only the implication case requires the distributivity assumption. 
This means that for arbitrary general intuitionistic team Kripke model, all $[\bot,\wedge,\vee]$-formulas are closed under finite joins. 

%
%

An immediate corollary  of the above lemma is that the idempotent law for $\vee$ and Axiom (2) are sound over distributive frames.

\begin{corollary}\label{gener_a2}
For any standard formula $\alpha$, and any formulas $\phi$ and $\psi$, we have that 
$\mathsf{DF}\models\alpha\vee\alpha\to\alpha$ and $\mathsf{DF}\models(\phi\to\alpha)\to((\psi\to\alpha)\to(\phi\vee\psi\to\alpha))$.
\end{corollary}
\begin{proof}
The soundness of $\alpha\vee\alpha\to\alpha$ follows from \Cref{Cup_closure_gts} and \Cref{eq_Ucl_idem}. The soundness of Axiom (2) follows from the observation that given Axiom (3), Axiom (2) is equivalent to $\alpha\vee\alpha\to\alpha$.
\end{proof}

%

As a consequence of Corollary \ref{distr_nec_cond} - Corollary \ref{gener_a2}, distributivity is a necessary and sufficient condition for validating Axiom (2), Axiom (3), idempotent law of $\vee$ for standard formulas and closure under join for standard formulas over finite classical frames. 

Finally, we show that over distributive frames the \textsf{Split} axiom is also sound, and thus all \TIPC axioms are recovered in $\mathsf{DF}$ or $\mathsf{DF}\models\TIPC$. We first prove the following lemma.

\begin{lemma}\label{gener_distr_G1}
Let $\mathfrak{F}=(W , R, \succcurlyeq, \Cup)$ be a general frame. For any teams  $t,r,s\subseteq W$, if $t\succcurlyeq\circ R r$ and $t\succcurlyeq \circ Rs$, then $R[t]\succcurlyeq r\Cup s$. 
\end{lemma}
\begin{proof}
By condition (b) we have that $R[t]\succcurlyeq r$ and $R[t]\succcurlyeq s$, which then imply 
that $R[t]\succcurlyeq r\Cup s$.
\end{proof}


\begin{proposition}\label{gener_split}
For any standard formula $\alpha$, and any formulas $\phi$ and $\psi$, we have that 
$\mathsf{DF}\models(\alpha\to\phi\vvee\psi)\to(\alpha\to\phi)\vvee(\alpha\to\psi)$. 
\end{proposition}
\begin{proof}
Let  $\MM=(W,R,\succcurlyeq,\Cup,V)$ be an arbitrary model with a distributive underlying frame. Suppose $\MM,t\not\models (\alpha\to\phi)\vvee(\alpha\to\psi)$. We show that $\MM,t\not\models\alpha\to\phi\vvee\psi$. By assumption, there are teams $r,s\subseteq W$ such that $t\succcurlyeq\circ R r$, $t\succcurlyeq\circ R s$, $r\models \alpha$, $r\not\models\phi$, $s\models\alpha$ and $s\not\models\psi$. Consider the team $r\Cup s$. By \Cref{Cup_closure_gts}, we have $r\Cup s\models\alpha$.
Since  $r,s\preccurlyeq r\Cup s$,  by persistence, we have $r\Cup s\not\models\phi$ and $r\Cup s\not\models\psi$, meaning $r\Cup s\not\models\phi\vvee\psi$. 
Now, by Lemma \ref{gener_distr_G1}, we have $R[t]\succcurlyeq r\Cup s$, from which we conclude that $R[t]\not\models\alpha\to\phi\vvee\psi$. Finally, since $tR^\circ R[t]$, we obtain $t\not\models\alpha\to\phi\vvee\psi$ by persistence.
\end{proof}

Hence, we have shown that $\mathsf{DF}\models\TIPC$. Whether the logic of $\mathsf{DF}$ is exactly \TIPC, or what the logic of $\mathsf{DF}$ is is left as future work.

\section{Concluding remarks and open problems}\label{sec:conclusion}

In this paper, we have studied intermediate logics in the setting of team semantics. Starting from the team-based intuitionistic logic \TIPC introduced in \cite{CiardelliIemhoffYang20}, we have explored two alternative approaches to define intermediate logics for team semantics. In the first approach, we extend \TIPC with axioms written in the standard language only. We showed that universal models for \IPC behave also as universal models for \TIPC, and a large class of intermediate standard axioms (including some de Jongh-formulas) still define the same class of Kripke frames. In the second approach, we made a first attempt to define team-based logics in which the \textsf{Split} axiom is not sound. This is done by generalizing the standard team semantics through  general intuitionistic team Kripke frames $\mathfrak{F}=(W,R,\succcurlyeq,\Cup)$, where $(\wp(W),\succcurlyeq,\Cup,\varnothing)$ is an abstraction of the underlying structure $(\wp(W),\supseteq,\cup,\varnothing)$ of teams in the standard team semantics.

A few  remarks on our choice of definitions in the second approach are in order, especially in comparison with two related approaches considered by Pun\u{c}och\'{a}\u{r} in \cite{Puncochar16,Puncochar17} (as sketched in the introduction section). Instead of modifying the underlying structure of teams to $(\wp(W),\succcurlyeq,\Cup,\varnothing)$, it also makes sense to consider the structure $(I,\supseteq,\cup,\varnothing)$, where $I\subseteq \wp(W)$ is a set of ``admissible'' teams, as it was done in  \cite{Puncochar16}. 
Another appealing alternative generalization is to  treat teams simply as elements of the domain $W$ in Definition \ref{g_frame_df}; that is, instead of ``lifting'' the order $R$ to the power set level, one ``flattens'' teams to the single world level. Or, as considered in  \cite{Puncochar17}, teams are viewed as primitive entities in an arbitrary bounded semilattice $(A,\Cup,0)$ (i.e., a team is an element in the domain $A$). This approach, however, does not seem to be in line with the fundamental idea of team semantics itself, at least not directly. Team semantics was introduced by Hodges \cite{Hodges1997a,Hodges1997b} to characterize notions of dependence, which, as Hodges observed, can only make sense in the presence of {\em multitudes} or {\em sets} of possibles worlds (i.e., teams), rather than in single worlds.  We thus choose to maintain this most distinguishing feature of team semantics, and define teams only as elements in the power set of the domain $W$. The benefit of this team semantics approach will become apparent when we investigate notions of dependence in the intermediate logic framework laid down in this paper, which we leave as future work.


There is another main difference between our approach and the approaches in \cite{Puncochar16,Puncochar17}. Our general intuitionistic team Kripke frames $\mathfrak{F}=(W,R,\succcurlyeq,\Cup)$ are, in a sense, always two-layered ---- the base layer consists of a usual intuitionistic Kripke frame $(W,R)$, and the structure $(\wp(W),\succcurlyeq,\Cup,\varnothing)$ belongs to the team layer. The ``intermediateness'' can come from variations in both layers as well as their interactions. In particular, our first approach (Section \ref{sec:standard_team}) can be viewed as keeping the team layer fixed to $(\wp(W),\supseteq,\cup,\varnothing)$, while varying the underlying Kripke frame $(W,R)$. The interaction between the two layers (as reflected in the composition $\succcurlyeq\circ R^\circ$) also affects the behavior of implication $\to$ (see Definition \ref{gen_team_sem_df}).  On the other hand, the approach in either \cite{Puncochar16} or \cite{Puncochar17} is essentially single-layered. The focus there is only on the team layer, being either $(I,\supseteq,\cup,\varnothing)$ or $(A,\succcurlyeq,\Cup,0)$ (where $\preccurlyeq$ is the partial order associated with the join $\Cup$), and the implication $\to$ is thus defined by using the partial order $\supseteq$ or $\succcurlyeq$ alone as the successor relation. In particular, putting $A=\wp(W)$, the setting in \cite{Puncochar17} essentially corresponds to the case $R=\textsf{id}$ in our setting. It is interesting to note, though, that in \cite{Puncochar17}, it was shown that over the frames $(A,\succcurlyeq,\Cup,0)$ with the semilattice $(A,\Cup,0)$ being distributive, the logic in the standard language (with local disjunction $\vee$) is exactly  \IPC. Distributivity is also shown to be a necessary and sufficient condition for join closure in the setting. Moreover, \cite{Puncochar17} also introduced a sound and complete system of natural deduction  for the language of \TIPC extended with the necessity $\Box$ modality, which is roughly \TIPC axioms/rules together with rules of \IK. These results correspond to our results in Section \ref{sec:recover}. Whether the arguments for the completeness theorem in \cite{Puncochar17} can be adapted in our context to show that distributivity is actually the characterization condition for \TIPC axioms  is open.  



A subtle but important technical difference between our setting and the setting in \cite{Puncochar17} also deserves comments. The semantics of the local disjunction $\vee$ in \cite{Puncochar17} is defined with the same clause as in our Definition \ref{gen_team_sem_df} except that our condition ``$t\preccurlyeq s\Cup r$'' is changed to ``$t=s\Cup r$" in \cite{Puncochar17}. While the two definitions for $\vee$ are actually equivalent in the standard team semantics setting (when $\succcurlyeq=\supseteq$ and $\Cup=\cup$, see Corollary \ref{local_v_alternative_df}), our definition turns out to be more general in the generalized team semantics setting. In particular, the persistence property  in \cite{Puncochar17} (especially for the $\vee$ case in the proof) also requires the assumption that the underlying semilattice in question is distributive (in which case \textsf{Split} axiom is also sound), whereas the persistence property in our setting (Lemma \ref{more_general_persistence}) does not require the distributivity assumption (and thus the \textsf{Split} axiom is in general not sound in our setting).




We end the paper by mentioning some more open problems. \Cref{sec:gt} is a first step towards a generalized team semantics and many further directions deserve investigation. For instance, in \Cref{g_frame_df} of general intuitionistic team Kripke frames, whether it is possible to weaken conditions (b) and (c) and still preserve the persistence property of the logic needs to be further explored. It also makes sense to weaken the assumption that $\preccurlyeq$ is the partial order of the join-semilattice $(\wp(W),\Cup,\varnothing)$. While this lattice order assumption was  used heavily in Section \ref{sec:recover} to recover all \TIPC axioms, in \Cref{sec:gtipl} the assumption is used actually only in establishing \Cref{gener_fact_vvee}(\ref{gener_fact_vvee_semilattice}) (the associativity of the local disjunction $\vee$), where the minimal condition seems to be that ``$t \preccurlyeq s$ implies $t\Cup r\preccurlyeq s\Cup r$ for all teams $r,s,t$''. 
Weakening the lattice order assumption to this alternative assumption will give rise to a different logic that requires further investigation. How to axiomatize the logic of the class of all general intuitionistic team Kripke frames,  what is the logic of the class $\mathsf{DF}$ of all distributive such frames, and do these logics have the disjunction property are also topics for future research. 
Providing appropriate definitions of universal models and de Jongh formulas in the generalized team semantics setting is also an interesting direction. Related work on universal models for fragments of intuitonistic modal logic
has been studied recently in \cite{BBCGGJ21}.
 Another very recent work \cite{Quadrellaro21} has provided algebraic semantics for \TIPC and the team-based intermediate logics obtained from our first approach (as defined in \Cref{team_intermediate_df1}). 
 Some other work on algebraic semantics of inquisitive and other relevant logics can be found in \cite{BezhanishviliGrillettiQuadrellaro20,GrillettiQuadrellaro21}. It is natural to study also the algebraic semantics for the logics obtained from our second approach. A recent Master's thesis \cite{DitrievaMsc} introduced a related type of general team semantics for a positive (modal) logic by considering complete (modal) lattices (where the distributivity of the lattices also plays a related subtle role). 
Studying modal extensions of the language \TIPC is also an interesting  direction for future work. Different variants of modal dependence and inquisitive logic have been studied by many authors, see, e.g., \cite{KontinenMulleerSchnoorVollmer15,sevenster09,VaMDL08,Yang17MD}.


%

\section*{Acknowledgments}

The authors are thankful to the referee for useful comments. They would also like to thank Davide Quadrellaro for interesting discussion related to the topic of this paper. 
The second author was supported by grant 330525 of the Academy of Finland, and Research Funds of the University of Helsinki.

\bibliographystyle{plain}

\end{document}